\documentclass[12pt,reqno]{amsart}
\usepackage{amssymb,amscd}   

\newtheorem{theorem}{Theorem}[section]
\newtheorem{proposition}[theorem]{Proposition} 
\newtheorem{corollary}[theorem]{Corollary}
\newtheorem{lemma}[theorem]{Lemma}
\newtheorem{remark}[theorem]{Remark}

\newtheorem{definition}{Definition}[section]

\def\Ps{\Psi}
\def\G{\Gamma}

\begin{document}
\title[Treeablilty and Costs for Discrete Groupoids]{Notes on Treeablilty and Costs for Discrete Groupoids in Operator Algebra Framework}
\author[Y. Ueda]
{Yoshimichi UEDA}
\thanks{Supported in part by Grant-in-Aid for Young Scientists (B) 14740118.}
\address{Graduate School of Mathematics, Kyushu University, 
Fukuoka, 810-8560, Japan}
\email{ueda@math.kyushu-u.ac.jp}
\maketitle

\section{Introduction} 

These notes discuss recent topics in orbit equivalence theory in operator algebra framework. Firstly, we provide an operator algebraic interpretation of discrete measurable groupoids in the course of giving a simple observation, which re-proves (and slightly generalizes) a result on treeability due to Adams and Spatzier \cite[Theorem 1.8]{Adams-Spatzier:AmerJMath90}, by using operator algebra techniques. Secondly, we reconstruct Gaboriau's work \cite{Gaboriau:InventMath00} on costs of equivalence relations in operator algebra framework with avoiding any measure theoretic argument. It is done in the same sprit as of \cite{Popa:MathScand85} for aiming to make Gaboriau's beautiful work much more accessible to operator algebraists (like us) who are not much familiar with ergodic theory. As simple byproducts, we clarify what kind of results in \cite{Gaboriau:InventMath00} can or cannot be generalized to the non-principal groupoid case, and observe that the cost of a countable discrete group with regarding it as a groupoid (i.e., a different quantity from Gaboriau's original one \cite[p.43]{Gaboriau:InventMath00}) is nothing less than the smallest number of its generators in sharp contrast with the corresponding $\ell^2$-Betti numbers, see Remark \ref{Rem3-4} (2). The methods given here may be useful for further discussing the attempts, due to Shlyakhtenko \cite{Shlyakhtenko:CMP01}\cite{Shlyakhtenko:Duke03}, of interpreting Gaboriau's work on costs by the idea of free entropy (dimension) due to Voiculescu. 

We introduce the notational convention we will employ; for a von Neumann algebra $N$, the unitaries, the partial isometries and the projections in $N$ are denoted by $N^u$, $N^{pi}$ and $N^p$, respectively. The left and right support projections of $v \in N^{pi}$ are denoted by $l(v)$ and $r(v)$, respectively, i.e., $l(v):=vv^*$ and $r(v):=v^* v$. We also mention that only von Neumann algebras with separable preduals will be discussed throughout these notes. 

We should thank Damien Gaboriau who earnestly explained us the core idea in his work, and also thank Tomohiro Hayashi for pointing out an insufficient point in a preliminary version. The present notes were  provided in part for the lectures we gave at University of Tokyo, in 2004, and we thank Yasuyuki Kawahigashi for his invitation and hospitality. 

\section{A Criterion for Treeability} 

Let $M \supseteq A$ be an inclusion of (not necessarily finite) von Neumann algebras with a faithful normal conditional expectation $E_A^M : M \rightarrow A$. Let $\mathbb{K}(M\supseteq A)$ be the $C^*$-algebra obtained as the operator norm $\Vert\ \cdot\ \Vert_{\infty}$-closure of $M e_A M$ on $L^2\left(M\right)$ with the Jones projection $e_A$ associated with $E_A^M$, and it is called the algebra of relative compact operators associated with the triple $M \supseteq A$, $E_A^M$. We use the notion of Relative Haagerup Property due to Boca \cite{Boca:Pacific93}. (Quite recently, Popa used  a slightly different formulation of Relative Haagerup Property in the type II$_1$ setting, see \cite{Popa:Preprint'03}, but we employ Boca's in these notes.) The triple $M \supseteq A$, $E_A^M$ is said to have Relative Haagerup Property if there is a net of $A$-bimodule (unital) normal completely positive maps $\Ps_{\lambda} : M \rightarrow M$, $\lambda \in \Lambda$, with $E_A^M\circ \Ps_{\lambda} = E_A^M$ for every $\lambda \in \Lambda$ such that for a fixed (and hence any) faithful state $\varphi \in M_*$ with $\varphi\circ E_A^M=\varphi$ one has  
\begin{itemize} 
\item $\lim_{\lambda}\left\Vert\Ps_{\lambda}(x) - x \right\Vert_{\varphi} = 0$ for every $x \in M$, or equivalently $\lim_{\lambda} \Ps_{\lambda} = \mathrm{id}_M$ pointwisely in $\sigma$-strong topology;  
\item $\widehat{\Ps}_{\lambda} \in \mathbb{K}(M\supseteq A)$,  
\end{itemize}
where $\widehat{\Ps}_{\lambda}$ is the bounded operator on $L^2\left(M\right)$ defined by $\widehat{\Ps}_{\lambda}\Lambda_{\varphi}(x) := \Lambda_{\varphi}\left(\Ps_{\lambda}(x)\right)$ for $x \in M$ with the canonical injection $\Lambda_{\varphi} : M \rightarrow L^2\left(M\right)$.  

The next lemma can be proved in the essentially same way as in \cite{Choda:ProcJapanAcad:83}, where group von Neumann algebras are dealt with. Although the detailed proof is now available in \cite[Proposition 3.5]{Jolissaint:JOT02}, we give its sketch for the reader's convenience, with focusing the ``only if" part, which we will need later. 

\begin{lemma}\label{Lem2-1} Assume that $M = A\rtimes_{\alpha}G$, i.e., $M$ is the crossed-product of $A$ by an action $\alpha$ of a countable discrete group $G$. Suppose that the action $\alpha$ has an invariant faithful state $\phi \in A_*$. Then, the inclusion $M \supseteq A$ with the canonical conditional expectation $E_A^M : M \rightarrow A$ has Relative  Haagerup Property if and only if $G$ has Haagerup Property {\rm (}see {\rm \cite{Haagerup:InventMath79},\cite{Choda:ProcJapanAcad:83})}. 
\end{lemma}

\begin{proof} (Sketch) Let $\lambda_g$, $g\in G$, be the canonical generators of $G$ in $M = A\rtimes_{\alpha}G$. The ``if" part is the easier implication. In fact, if $G$ has Haagerup Property, i.e., there is a net of positive definite functions $\psi_{\lambda}$ with vanishing at infinity such that $\psi_{\lambda}(g) \rightarrow 1$ for every $g \in G$, then the required $\Ps_{\lambda}$ can be constructed in such a way that $\Ps_{\lambda}\big(\sum_g a(g)\lambda_g\big) := \sum_{g\in G} \psi_{\lambda}(g)a(g)\lambda_g$ for every finite linear combination $\sum_{g\in G} a(g)\lambda_g \in A\rtimes_{\alpha}G$, see \cite[Lemma 1.1]{Haagerup:InventMath79}. The ``only if" part is as follows. Define $\psi_{\lambda}(g) := \phi\circ E_A^M\big(\Ps_{\lambda}(\lambda_g)\lambda_g^*\big)$, $g \in G$, and clearly $\psi_{\lambda}(g) \rightarrow 1$ for every $g \in G$. Let $\varepsilon>0$ be arbitrary small.  One can choose a $T = \sum_{i=1}^n x_i e_A y_i \in Me_A M$ with $\big\Vert \widehat{\Ps}_{\lambda} - T\big\Vert_{\infty}\leq\varepsilon/2$. Then, for $g \in G$ one has 
$\big|\psi_{\lambda}(g)\big| \leq \frac{\varepsilon}{2} + 
\sum_{i=1}^n \big\Vert x_i \big\Vert_{\infty} \big\Vert E_A^M(y_i\lambda_g)\big\Vert_{\phi}$.
Since $\big\Vert y_i \big\Vert_{\phi\circ E_A^M}^2 = \sum_{h \in G}\left\Vert E_A^M\left(y_i\lambda_h^*\right)\right\Vert_{\phi}^2$ (where it is crucial that $\phi$ is invariant under $\alpha$), one can choose a finite subset $K$ of $G$ in such a way that every $g \in G\setminus K$ satisfies that $\big\Vert E_A^M\big(y_i\lambda_g\big)\big\Vert_{\phi}\leq\varepsilon/\big(2\sum_{i=1}^n \Vert x_i \Vert_{\infty}\big)$. Then, $|\psi_{\lambda}(g)|\leq\varepsilon$ for every $g \in G\setminus K$. 
\end{proof} 
  
In what follows, we further assume that $A$ is commutative. Denote $\mathcal{G}(M\supseteq A) := \big\{ v \in M : v^*v, vv^* \in A^p,\ vAv^*=Avv^* \big\}$ and call it  the full (normalizing) groupoid of $A$ in $M$. When $A$ is a MASA in $M$ and $\mathcal{G}(M\supseteq A)$ generates $M$ as von Neumann algebra, we call $A$ a Cartan subalgebra in $M$, see \cite{Feldman-Moore:TransAMS77}. Let us introduce a von Neumann algebraic formulation of the set of one-sheeted sets in a countable discrete measurable groupoid. 
\begin{definition}\label{Def2-1} An $E_A^M$-groupoid is a subset $\mathcal{G}$ of $\mathcal{G}(M\supseteq A)$ equipped with the following properties: 
\begin{itemize} 
\item $u, v \in \mathcal{G} \Longrightarrow uv \in \mathcal{G}${\rm ;}  
\item $u \in \mathcal{G} \Longrightarrow u^* \in \mathcal{G}${\rm ;}
\item $u \in A^{pi} \Longrightarrow u \in \mathcal{G}$ {\rm (}and, in particular, $u \in \mathcal{G}, p \in A^p \Longrightarrow pu, up \in \mathcal{G}${\rm );}
\item Let $\left\{u_k\right\}_k$ be a {\rm (}possibly infinite{\rm )} collection of elements in $\mathcal{G}$. If the support projections and the range projections respectively form mutually orthogonal families, then $\sum_k u_k \in \mathcal{G}$ in $\sigma$-strong* topology{\rm ;} 
\item Each $u \in \mathcal{G}$ has a {\rm (}possibly zero{\rm )} $e \in A^p$ such that $e \leq l(u)$ and  $E_A^M(u)=eu${\rm ;}
\item Each $u \in \mathcal{G}$ satisfies that $E_A^M(uxu^*) = uE_A^M(x)u^*$ for every $x \in M$. 
\end{itemize} 
\end{definition} 

Such a projection $e$ as in the fifth is uniquely determined as the modulus part of the polar decomposition of $E_A^M(u)$. The sixth automatically holds, either when $E_A^M$ is the (unique) $\tau$-conditional expectation with a faithful tracial state $\tau \in M_*$ or when $A$ is a MASA in $M$. The next two lemmas are proved based on the same idea as for \cite[Proposition 2.2]{Popa:MathScand85}.  

\begin{lemma}\label{Lem2-2} Let $\mathcal{U}$ be an {\rm (}at most countably infinite{\rm )} collection of elements in $\mathcal{G}$, and $w_0 :=1, w_1,\dots$ be the words in $\mathcal{U}\sqcup\mathcal{U}^*$ of reduced form in the formal sense with regarding $u^{-1}= u^*$ for $u \in \mathcal{U}$. Suppose that  $\mathcal{G}'' = A\vee\mathcal{U}''$ as von Neumann algebra. Then, each $v \in \mathcal{G}$ has a partition $l(v) = \sum_k p_k$ in $A^p$ with $p_k v= p_k E_A^M\left(vw_k^*\right)w_k$ for every $k$. Furthermore, each coefficient $E_A^M\left(vw_k^*\right)$ falls in $A^{pi}$ and $v = \sum_k p_k E_A^M\left(vw_k^*\right) w_k$ in $\sigma$-strong* topology. 
\end{lemma}    
\begin{proof} By the fifth requirement of $E_A^M$-groupoids one can find a (unique) $e_k \in A^p$ in such a way that $e_k \leq l(v)$ and $E_A^M\left(vw_k^*\right) = e_k vw_k^*$, i.e., $e_k v = E_A^M\left(vw_k^*\right)w_k$. Set $e := \bigvee_k e_k$, and choose a faithful state $\varphi \in M_*$ with $\varphi\circ E_A^M = \varphi$. Then, we get $\big(v -ev| aw_k\big)_{\varphi}:=\varphi\big((aw_k)^*(v-ev)\big)= 0$ for every $a \in A$ and every $k$, where the sixth requirement is used crucially. Since the $aw_k$'s give a total subset in  $\mathcal{G}''$ in $\sigma$-strong topology, we concludes that $v = ev$ so that $e=l(v)$. Since $A$ is commutative, one can construct $p_0, p_1,\dots \in A^p$ in such a way that $p_k \leq e_k$ and $\sum_k p_k = e$. Then, we have $p_k v = p_k e_k v = p_k E_A^M\left(vw_k^*\right)w_k$ for each $k$. \end{proof}  

\begin{lemma}\label{Lem2-2.1} Let $\mathcal{G}_0 \subseteq \mathcal{G}$ be an $E_A^M$-groupoid with a faithful normal conditional expectation $E_{M_0}^M : M \rightarrow M_0 := \mathcal{G}_0''$ with $E_A^M\circ E_{M_0}^M = E_A^M$. Then, each $u \in \mathcal{G}$ has a {\rm (}unique{\rm )} $p \in A^p$ such that $p \leq l(u)$ and $E_{M_0}^M(u)=pu$. 
\end{lemma}   
\begin{proof} By the same method as for the previous lemma together with standard exhaustion argument one can construct an (at most countably infinite) subset $\mathcal{W}$ of $\mathcal{G}_0$ that possesses the following properties: $E_A^M\left(w_1 w_2^*\right) = \delta_{w_1,w_2} l(w_1)$ for $w_1, w_2 \in \mathcal{W}$; each $u \in \mathcal{G}$ has an orthogonal family $\{e_w(u)\}_{w\in\mathcal{W}}$ such that $e_w(u) \leq l(uw^*)$ ($\leq l(u)\wedge r(w)$) and $e_w(u)u= E_A^M(uw^*)w$; if $u \in \mathcal{G}$ is chosen from $\mathcal{G}_0$, then $l(u) = \sum_{w\in\mathcal{W}} e_w(u)$. Choose a faithful state $\varphi \in M_*$ with $\varphi\circ E_A^M = \varphi$. Set $p := \sum_{w\in\mathcal{W}} e_w(u) \leq l(u)$, and we have $\big(E_{M_0}^M(u)-pu\big|aw\big)_{\varphi}:= \varphi\big((aw)^*(E_{M_0}^M(u)-pu)\big)= 0$ for every $aw$, $a \in A, w \in \mathcal{W}$. Hence, we get $E_{M_0}^M(u)=pu$. The uniqueness follows from that for the (right) polar decomposition of $E_{M_0}^M(u)$.     
\end{proof} 

Let $\begin{matrix} \mathcal{G}_{11} & \supset & \mathcal{G}_{12} \\ 
\cup & & \cup \\ \mathcal{G}_{21} & \supset & \mathcal{G}_{22}
\end{matrix}$ be $E_A^M$-groupoids and write $M_{ij} := \mathcal{G}_{ij}''$. Assume that there are faithful normal conditional expectations $E_{ij} : M \rightarrow M_{ij}$ with $E_A^M\circ E_{ij} = E_A^M$. In this case, Lemma \ref{Lem2-2.1} enables us to see that the following three conditions are equivalent:  $\begin{matrix} M_{11} & \supset & M_{12} \\ 
\cup & & \cup \\ M_{21} & \supset & M_{22} \end{matrix}$ forms a commuting square; $M_{22} = M_{12}\cap M_{21}$; and $\mathcal{G}_{22} = \mathcal{G}_{12}\cap\mathcal{G}_{21}$. Moreover, one also observes, in the similar way as above, that if two $E_A^M$-groupoids inside a fixed $\mathcal{G}$ generate the same intermediate von Neumann algebra between $\mathcal{G}'' \supseteq A$, then they must coincide. If $A = \mathbf{C}1$, then the image $\pi(\mathcal{G})$ with the quotient map $\pi : M^u \rightarrow M^u/\mathbb{T}1$ is a countable discrete group. The full groupoid $\mathcal{G}(M\supseteq A)$ itself becomes an $E_A^M$-groupoid when $A$ is a MASA in $M$ thanks to Dye's lemma (\cite[Lemma 6.1]{Dye:AmerJMath63}; also see \cite{Choda:ProcJapanAcad:65}), which asserts the same as in Lemma \ref{Lem2-2.1} for $\mathcal{G}(M\supseteq A)$ without any assumption when $A$ is a Cartan subalgebra in $M$. (The non-finite case needs a recently well-established result in \cite{Aoi:JMSJ03}.) Moreover, the set of one-sheeted sets in a countable discrete measurable groupoid canonically gives an $E_A^M$-groupoid, where $M \supseteq A$ with $E_A^M : M \rightarrow A$ are constructed by the so-called regular representation. See just after the next lemma for this fact. Let us introduce the notions of graphings and  treeings due to Adams \cite{Adams:ErgodicTheoryDynamSystems:90} (also see \cite{Gaboriau:InventMath00}, \cite[Proposition 7.5]{Shlyakhtenko:CMP01}) in operator algebra framework. We call such a collection $\mathcal{U}$ as in Lemma \ref{Lem2-2}, i.e, $\mathcal{G}''=A\vee\mathcal{U}''$, a graphing of $\mathcal{G}$. On the other hand, a collection $\mathcal{U}$ of elements in $\mathcal{G}(M\supseteq A)$ ({\it n.b.}, not assumed to be a graphing) is said to be a treeing if $E_A^M(w)=0$ for all words $w$ in $\mathcal{U}\sqcup\mathcal{U}^*$ of reduced form in the formal sense. This is equivalent to that $\mathcal{U}$ is a $*$-free family (or equivalently, $\left\{ A\vee\{u\}''\right\}_{u\in\mathcal{U}}$ is a free family of von Neumann algebras) with respect to $E_A^M$ in the sense of Voiculescu (see e.g.~\cite[\S\S3.8]{VDN:Book92}) since every element in $\mathcal{G}(M\supseteq A)$ normalizes $A$. We say that $\mathcal{G}$ has a treeing $\mathcal{U}$ when $\mathcal{U}$ is a treeing and a graphing of $\mathcal{G}$, and also $\mathcal{G}$ is treeable if $\mathcal{G}$ has a treeing. 

\begin{lemma}\label{Lem2-3} {\rm (cf.~\cite{Haagerup:InventMath79}, \cite{Boca:Pacific93})} If an $E_A^M$-groupoid $\mathcal{G}$ has a treeing $\mathcal{U}$, then the inclusion $M(\mathcal{G}) := \mathcal{G}'' \supseteq A$ with $E_A^M\big|_{M(\mathcal{G})} : M(\mathcal{G}) \rightarrow A$ must have Relative Haagerup Property. 
\end{lemma}  
\begin{proof} We may and do assume $M=M(\mathcal{G})$ for simplicity. We first assume that $\mathcal{U}$ is a finite collection. Since $\mathcal{U}$ is a treeing, each $u \in \mathcal{U}$ satisfies either $u^m\neq0$ or $E_A^M\left(u^m\right)=0$, and thus each $N_u := A\vee\{u\}''$ can be decomposed into    
\begin{equation*} 
\text{(i)}\quad N_u = \sideset{}{^{\oplus}}\sum_{|m|\leq n_u} u^m A \quad \text{or} \quad 
\text{(ii)}\quad N_u = \sideset{}{^{\oplus}}\sum_{m\in\mathbb{Z}} u^m A
\end{equation*} 
in the Hilbert space $L^2(M)$ via $\Lambda_{\varphi}$ with a faithful state $\varphi \in M_*$ with $\varphi\circ E_A^M=\varphi$. Here, $u^{-m}$ means the adjoint $u^*{}^m$ as convention. By looking at this description, it is not so hard to confirm that each triple $N_u \supseteq A$ with $E_A^M\big|_{N_v}$ satisfies Relative Haagerup Property. Namely, one can construct a net $\Ps_u^{(\varepsilon)} : N_u \rightarrow N_u$ of completely positive maps in such a way that 
\begin{itemize} 
\item $E_A^M\circ \Ps_u^{(\varepsilon)} = E_A^M\big|_{N_u}$; 
\item $\Ps_u^{(\varepsilon)}$ converges to $\mathrm{id}_{N_u}$ pointwisely, in 
$\sigma$-strong topology, as $\varepsilon \searrow 0$; 
\item $\widehat{\Ps}_u^{(\varepsilon)}$ falls into $\mathbb{K}\left(N_u\supseteq A\right)$ in $L^2(N_u)= \overline{\Lambda_{\varphi}(N_u)}$; 
\item $T_u^{(\varepsilon)} := \widehat{\Ps}_u^{(\varepsilon)}\big|_{L^2(N_v)^{\circ}}$ satisfies $\left\Vert T_u^{(\varepsilon)}\right\Vert_{\infty} = \exp(-\varepsilon)$ with $L^2(N_u)^{\circ} := (1-e_A)L^2(N_u)$.  
\end{itemize} 
The case (i) is easy, that is, 
\begin{equation*} 
\Ps_u^{(\varepsilon)} := e^{-\varepsilon}\ \mathrm{id}_{N_u} + (1-e^{-\varepsilon}) E_A^M\big|_{N_u} = E_A^M\big|_{N_u} + e^{-\varepsilon}\left(\mathrm{id}_{N_u}-E_A^M\big|_{N_u}\right)  
\end{equation*} 
converges to $\mathrm{id}_{N_u}$ pointwidely, in $\sigma$-strong topology, and one has 
\begin{equation}\label{eq1-L2-3} 
\widehat{\Ps}_u^{(\varepsilon)} = e_A + e^{-\varepsilon}\left(\sum_{0\lneqq|m|\leq n_u} u^m e_A u^{-m}\right) \in N_u e_A N_u. 
\end{equation} 
The case (ii) needs to modify the standard argument \cite[Lemma 1.1]{Haagerup:InventMath79}. By using the cyclic representation of $\mathbb{Z}$ induced by the positive definite function $m \mapsto e^{-\varepsilon|m|}$ one can construct a sequence $s_k \in \ell^{\infty}(\mathbb{Z})$ satisfying that $\sum_k \left|s_k(m)\right|^2 < +\infty$ for every $m \in \mathbb{Z}$ and moreover that $\sum_k s_k(m_1)\overline{s_k(m_2)} = e^{-\varepsilon(|m_1-m_2|)}$ for every pair $m_1, m_2 \in \mathbb{Z}$. Set $S_k := \sum_{m\in\mathbb{Z}} s_k(m) u^m e_A u^{-m}$ (on $L^2(N_u)$), and then the desired completely positive maps can be given by 
\begin{equation*} 
\Ps_u^{(\varepsilon)} : x \in N_u \mapsto \sum_k S_k x S_k^* \in B\left(L^2\left(N_u\right)\right). 
\end{equation*}
(Note here that $u^m e_A u^{-m}$ is the projection from $L^2(N_u)$ onto $\overline{\Lambda_{\varphi}\left(u^m A\right)}$.) In fact, it is easy to see that $\Ps_u^{(\varepsilon)}\left(u^m a\right) = e^{-\varepsilon|m|} u^m a$ for $m \in \mathbb{Z}$, $a \in A$, which shows that the range of $\Ps_u^{(\varepsilon)}$ sits in $N_u$ and that $\Ps_u^{(\varepsilon)}$ converges to $\mathrm{id}_{N_u}$ pointwidely, in $\sigma$-strong topology. Moreover, one has  
\begin{equation}\label{eq2-L2-3}  
\widehat{\Ps}_u^{(\varepsilon)} = \sum_{m\in\mathbb{Z}} e^{-\varepsilon|m|} u^m e_A u^{-m} 
= \lim_{n\rightarrow\infty} \sum_{|m|\leq n} e^{-\varepsilon|m|} u^m e_A u^{-m}   
\end{equation}   
in operator norm.  

Since $\mathcal{U}$ is a treeing, we have 
\begin{equation*}
\left(M, E_A^M\right) = \underset{u\in\mathcal{U}}{\bigstar_A}\left(N_u, E_A^M\big|_{N_u}\right). 
\end{equation*}
Therefore, \cite[Proposition 3.9]{Boca:Pacific93} shows that the inclusion  $M \supseteq A$ with $E_A^M$ satisfies Relative Haagerup Property since we have shown that so does each $N_u \supseteq A$ with $E_A^M\big|_{N_u}$. However, we would like to give the detailed argument on this point for the reader's convenience. Thanks to $E_A^M\circ \Ps_u^{(\varepsilon)} = E_A^M\big|_{N_u}$, we can construct the free products of completely positive maps $\Ps^{(\varepsilon)} := \underset{u\in\mathcal{U}}{\bigstar_A} \Ps_u^{(\varepsilon)} : M \rightarrow M$, which is uniquely determined by the following properties: 
\begin{itemize} 
\item $E_A^M\circ\Ps^{(\varepsilon)} = E_A^M$; 
\item $\Ps^{(\varepsilon)}\left(x_1 x_2 \cdots x_{\ell}\right) = 
\Ps_{u_1}^{(\varepsilon)}\left(x_1\right)\Ps_{u_2}^{(\varepsilon)}\left(x_2\right)\cdots\Ps_{u_{\ell}}^{(\varepsilon)}\left(x_{\ell}\right)$ for $x_j^{\circ} \in \mathrm{Ker}E_A^M\cap N_{u_j}$ with $u_1\neq u_2\neq\cdots\neq u_{\ell}$. 
\end{itemize}  
(See \cite[Theorem 3.8]{BlanchardDykema:Pacific01} in the most generic form at present.) Since each $\Ps_u^{(\varepsilon)}$ converges to $\mathrm{id}_{N_u}$ pointwidely in $\sigma$-strong topology, as $\varepsilon\searrow0$, the above two properties enable us to confirm that so does $\Ps^{(\varepsilon)}$ to $\mathrm{id}_M$. It is standard to see that 
\begin{equation*} 
\widehat{\Ps}^{(\varepsilon)} = 1_{L^2(A)}\oplus\sideset{}{^{\oplus}}\sum_{\ell\geq1}\sideset{}{^{\oplus}}\sum_{u_1\neq u_2 \neq \cdots \neq u_{\ell}}T_{u_1}^{(\varepsilon)}\otimes_{\varphi}T_{u_2}^{(\varepsilon)}\otimes_{\varphi}\cdots\otimes_{\varphi}T_{u_{\ell}}^{(\varepsilon)}
\end{equation*} 
in the free product representation 
\begin{equation}\label{eq3-L2-3} 
L^2(M) = L^2(A)\oplus\sideset{}{^{\oplus}}\sum_{\ell\geq1}\sideset{}{^{\oplus}}\sum_{u_1\neq u_2 \neq \cdots \neq u_{\ell}}L^2\left(N_{u_1}\right)^{\circ}\otimes_{\varphi}\cdots\otimes_{\varphi} L^2\left(N_{u_{\ell}}\right)^{\circ} 
\end{equation} 
with $L^2(A) = \overline{\Lambda_{\varphi}(A)} \subseteq L^2(M)$, where $\otimes_{\varphi}$ means the relative tensor product operation over $A$ with respect to $\varphi|_A \in A_*$ (see \cite{Sauvageot:JOT83}). Notice that, with $x_j^{\circ} \in \mathrm{Ker}E_A^M\cap N_{u_j}$, $u_1\neq u_2\neq\cdots\neq u_{\ell}$,   
\begin{equation*}
\Lambda_{\varphi}\left(x_1^{\circ} x_2^{\circ} \cdots x_{\ell}^{\circ}\right)=\Lambda_{\varphi}\left(x_1^{\circ}\right)\otimes_{\varphi}\Lambda_{\varphi}\left(x_2^{\circ}\right)\otimes_{\varphi}\cdots\otimes_{\varphi}\Lambda_{\varphi}\left(x_{\ell}^{\circ}\right)
\end{equation*} 
in \eqref{eq3-L2-3}, and hence by \eqref{eq1-L2-3},\eqref{eq2-L2-3}, we have, via \eqref{eq3-L2-3},   
\allowdisplaybreaks{
\begin{align*} 
&\widehat{\Ps}^{(\varepsilon)}\big|_{L^2\left(N_{u_1}\right)^{\circ}\otimes_{\varphi}L^2\left(N_{u_2}\right)^{\circ}\otimes_{\varphi}\cdots\otimes_{\varphi} L^2\left(N_{u_{\ell}}\right)^{\circ}} \\
&=
T_{u_1}^{(\varepsilon)}\otimes_{\varphi}T_{u_2}^{(\varepsilon)}\otimes_{\varphi}\cdots\otimes_{\varphi}T_{u_{\ell}}^{(\varepsilon)} \\
&= 
\sum_{m_1,m_2,\dots,m_{\ell}} e^{-\varepsilon n}\ 
u_1^{m_1} u_2^{m_2}\cdots u_{\ell}^{m_{\ell}} e_A u_{\ell}^{-m_{\ell}}\cdots u_2^{-m_2} u_1^{-m_1}
\end{align*}
}with certain natural numbers $n=n(u_1,u_2,\dots,u_{\ell}; m_1,m_2,\dots,m_{\ell})$ that converges to  $+\infty$ as $|m_1|,|m_2|\dots,|m_{\ell}| \rightarrow \infty$ (as long as when it is possible to do so). Note also that
\allowdisplaybreaks{
\begin{align*} 
\left\Vert T_{u_1}^{(\varepsilon)}\otimes_{\varphi}T_{u_2}^{(\varepsilon)}\otimes_{\varphi}\cdots\otimes_{\varphi}T_{u_{\ell}}^{(\varepsilon)}\right\Vert_{\infty} &\leq 
\left\Vert T_{u_1}^{(\varepsilon)} \right\Vert_{\infty}\cdot\left\Vert T_{u_2}^{(\varepsilon)} \right\Vert_{\infty}\cdots\left\Vert T_{u_{\ell}}^{(\varepsilon)} \right\Vert_{\infty} \\
&= 
e^{-\ell\varepsilon} \longrightarrow 0 \quad \text{(as $\ell \rightarrow \infty$)}. 
\end{align*} 
}By these facts, $\widehat{\Ps}^{(\varepsilon)}$ clearly falls in the operator norm closure of $M e_A M$ since $\mathcal{U}$ is a finite collection. Hence, the net $\Ps^{(\varepsilon)}$ of completely positive maps on $M$ provides a desired one showing that the inclusion $M\supseteq A$ with the $E_A^M$ has Relative Haagerup Property. 

Next, we deal with the case that $\mathcal{U}$ is an infinite collection. In this case, one should at first choose a filtration $\mathcal{U}_1 \subseteq \mathcal{U}_2 \subseteq \cdots \nearrow \mathcal{U} = \bigcup_k \mathcal{U}_k$ by finite sub-collections. Then, instead of the above $\Ps^{(\varepsilon)}$ we consider the completely positive maps 
\begin{equation*} 
\Ps^{(\varepsilon)}_k := \left(\underset{u \in \mathcal{U}_k}{\bigstar_A}\Ps_u^{(\varepsilon)}\right)\circ E_{M_k}^M : M \rightarrow M_k := \bigvee_{u\in\mathcal{U}_k} N_u \left(=\underset{u\in\mathcal{U}_k}{\bigstar_A} N_u\right) \rightarrow M_k \subseteq M
\end{equation*} 
with $\varphi\circ E_A^M$-conditional expectation $E_{M_k}^M : M \rightarrow M_k$. Since $M_1 \subseteq M_2 \subseteq \cdots \nearrow M = \bigvee_k M_k$, the non-commutative Martingale convergence theorem \cite[Lemma 2]{Connes:JFA75} says that $E_{M_k}^M$ converges to  $\mathrm{id}_M$ pointwidely, in $\sigma$-strong topology, as $k\rightarrow\infty$, and so does  $\Ps^{(\varepsilon)}_k$ to $\mathrm{id}_M$ too, as $\varepsilon\searrow0$, $k\rightarrow\infty$. We easily see that  
 \begin{equation}\label{eq4-L2-3}  
\widehat{\Ps}^{(\varepsilon)}_k = 1_{L^2(A)}\oplus\sideset{}{^{\oplus}}\sum_{\ell\geq1}\sideset{}{^{\oplus}}\sum_{u_1\neq u_2 \neq \cdots \neq u_{\ell} \atop u_j \in \mathcal{U}_k}T_{u_1}^{(\varepsilon)}\otimes_{\varphi}T_{u_2}^{(\varepsilon)}\otimes_{\varphi}\cdots\otimes_{\varphi}T_{u_{\ell}}^{(\varepsilon)}
\end{equation} 
in \eqref{eq3-L2-3}. Note that the summation of each $\ell$th direct summand of \eqref{eq4-L2-3} is taken over the alternating words in the fixed finite collection $\mathcal{U}_k$ of length $\ell$, and thus the previous argument works for showing that $\widehat{\Ps}^{(\varepsilon)}_k$ falls into $\mathbb{K}(M\supseteq A)$. Hence, we are done. 
\end{proof} 

Here, we briefly summarize some basic facts on von Neumann algebras associated with countable discrete measurable groupoids, see e.g.~\cite{Hahn:TransAMS78},\cite{Renault:LNM793}. Let $\G$ be a countable discrete measurable groupoid with unit space $X$, where $X$ is a standard Borel space with a regular Borel measure. With a non-singular measure on $X$ under $\G$ one can construct, in a canonical way, a pair $M(\G)\supseteq A(\G)$ of von Neumann algebra and distinguished commutative von Neumann subalgebra with $A(\G) = L^{\infty}(X)$ and a faithful normal conditional expectation $E_{\G} : M(\G) \rightarrow A(\G)$, by the so-called regular representation of $\G$ due to Hahn \cite{Hahn:TransAMS78} (also see \cite[Chap.~II]{Renault:LNM793}), which generalizes Feldman-Moore's construction \cite{Feldman-Moore:TransAMS77} for countable discrete measurable equivalence relations. Denote by $\mathcal{G}_{\G}$ of $\G$ the set of ``one-sheeted sets in $\G$" or called ``$\G$-sets", i.e., measurable subsets of $\G$, on which the mappings $\gamma \in \G \mapsto \gamma\gamma^{-1}, \gamma^{-1}\gamma \in X$ are both injective. Note that $\mathcal{G}_{\G}$ becomes an inverse semigroup with product $E_1 E_2 := \{\gamma_1\gamma_2 : \gamma_1 \in E_1, \gamma_2 \in E_2, \gamma_1^{-1}\gamma_1 = \gamma_2\gamma_2^{-1} \}$ and inverse $E \mapsto E^{-1} := \{\gamma^{-1} : \gamma \in E\}$. Each $E \in \mathcal{G}_{\G}$ gives an element $u(E) \in \mathcal{G}\left(M(\G)\supseteq A(\G)\right)$ with the properties: Its left and right support projections $l\left(u(E)\right), r\left(u(E)\right)$ coincide with the characteristic functions on $EE^{-1}=\{\gamma\gamma^{-1} : \gamma \in E\}, E^{-1}E=\{ \gamma^{-1}\gamma : \gamma \in E \}$, respectively, in $L^{\infty}(X)$; The mapping $u : E \in \mathcal{G}_\G \mapsto u(E) \in \mathcal{G}\left(M(\G)\supseteq A(\G)\right)$ is an inverse semigroup homomorphism (being injective modulo null sets), where $\mathcal{G}\left(M(\G)\supseteq A(\G)\right)$ is equipped with the inverse operation $u \mapsto u^*$; $E_{\G}(u(E)) = eu(E)$ with the projection $e$ given by the characteristic function on $X\cap E$; $E_{\G}\big(u(E)xu(E)^*\big)=u(E)E_{\G}(x)u(E)^*$ for every $x\in M(\G)$. It is not difficult to see that $\mathcal{G}(\G) := A(\G)^{pi}u\left(\mathcal{G}_{\G}\right)= \left\{ au(E) \in \mathcal{G}\left(M(\G)\supseteq A(\G)\right) : a \in A(\G)^{pi}, E \in \mathcal{G}_{\G} \right\}$ is an $E_{\G}$-groupoid, which generates  $M(\G)$ as von Neumann algebra. An (at most countably infinite) collection $\mathcal{E}$ of elements in $\mathcal{G}_{\G}$ is called a graphing of $\G$ if it generates $\G$ as groupoid, or equivalently the smallest groupoid that contains $\mathcal{E}$ becomes $\G$. If no word in $\mathcal{E}\sqcup\mathcal{E}^{-1}$ of reduced form in the formal sense intersects with the unit space $X$ of strictly positive measure, then we call $\mathcal{E}$ a treeing of $\G$. Then, it is not hard to see the following two facts: (i) the collection $u(\mathcal{E})$ of $u(E) \in \mathcal{G}\left(M(\G)\supseteq A(\G)\right)$ with $E \in \mathcal{E}$ is a graphing of $\mathcal{G}({\G})$ if and only if $\mathcal{E}$ is a graphing of $\G$; and similarly, (ii) the collection $u(\mathcal{E})$ is a treeing of $\mathcal{G}({\G})$ if and only if $\mathcal{E}$ is a treeing of $\G$. With these considerations, the previous two lemmas immediately imply the following criterion for treeability:  

\begin{proposition}\label{Prop2-4} Relative Haagerup Property of $M(\G)\supseteq A(\G)$ with $E_{\G}$ is necessary for treeability of countable discrete measurable groupoid $\G$. In particular, any countably infinite discrete group without Haagerup Property has no treeable free action with finite invariant measure.  
\end{proposition}  

Note that this follows from a much deeper result due to Hjorth (see \cite[\S28]{Kechris-Miller:LNM04}) with the aid of Lemma \ref{Lem2-1} if a given $\G$ is principal or an equivalence relation. The above proposition clearly implies the following result of Adams and Spatzier:  

\begin{corollary}\label{Thm2-5} {\rm (\cite[Theorem 1.8]{Adams-Spatzier:AmerJMath90})} Any countably infinite discrete group of Property T admits no treeable free ergodic action with finite invariant measure. 
\end{corollary}   

\begin{remark}\label{Rem2-6} Note that the finite measure preserving assumption is very important in the above assertions. In fact, any countably infinite discrete group of Property T has an amenable free ergodic action without invariant finite measure {\rm (}e.g.~the boundary actions of some word-hyperbolic groups and the translation actions of discrete groups on themselves{\rm )}. 
\end{remark} 

\section{Operator Algebra Approach to Gaboriau's Results}  

We explain how to re-prove Gaboriau's results \cite{Gaboriau:InventMath00} on costs of equivalence relations (and slightly generalize them to the groupoid setting) in operator algebra framework, avoiding any measure theoretic argument. Throughout this section, we keep and employ the terminologies in the previous section. 

Let $\mathcal{E}$ be a graphing of a countable discrete measurable groupoid $\G$ with a non-singular probability measure $\mu$ on the unit space $X$. Following Levitt \cite{Levitt:ErgodicTheoryDynamSystems95} and Gaboriau \cite{Gaboriau:InventMath00} the $\mu$-cost of $\mathcal{E}$ is defined to be 
\begin{equation*} 
C_{\mu}(\mathcal{E}) := \sum_{E \in \mathcal{E}} \frac{\mu\left(EE^{-1}\right)+\mu\left(E^{-1}E\right)}{2}, 
\end{equation*}
and the $\mu$-cost of $\G$ by taking the infimum all over the graphings, that is, 
\begin{equation*} 
C_{\mu}(\G) := \inf \left\{ C_{\mu}(\mathcal{E}) : \text{$\mathcal{E}$ graphing of $\G$}\right\}.
\end{equation*}     
In fact, if $\G$ is a principal one (or equivalently a countable discrete equivalence relation) with an invariant probability measure $\mu$, the $\mu$-cost of graphings and that of $\G$ coincide with Levitt and Gaboriau's ones. 

Let $M\supseteq A$ be a von Neumann algebra and a distinguished commutative von Neumann subalgebra with a faithful normal conditional expectation $E_A^M : M \rightarrow A$, and $\mathcal{G}$ be an $E_A^M$-groupoid. For a faithful state $\varphi \in M_*$ with $\varphi\circ E_A^M = \varphi$, the $\varphi$-cost of a graphing $\mathcal{U}$ of $\mathcal{G}$ is defined to be 
\begin{equation*} 
C_{\varphi}(\mathcal{U}) := \sum_{u \in \mathcal{U}} \frac{\varphi\left(l(u)+r(u)\right)}{2}, 
\end{equation*} 
and that of $\mathcal{G}$ by taking the infimum all over the graphings of $\mathcal{G}$, that is, 
\begin{equation*} 
C_{\varphi}(\mathcal{G}) := \inf\left\{ C_{\varphi}(\mathcal{U}) : \text{$\mathcal{U}$ graphing of $\mathcal{G}$}\right\}. 
\end{equation*} 
We sometimes consider those cost functions $C_{\varphi}$ for both graphings and $E_A^M$-groupoids with the same equations even when $\varphi$ is not a state (but still normal and positive). When $\mathcal{G}=\mathcal{G}(\G)$, i.e., the canonical $E_{\G}$-groupoid associated with a countable discrete measurable groupoid $\G$, it is plain to verify that $C_{\varphi}(\mathcal{G}(\G)) = C_{\mu}(\G)$ with the state $\varphi \in M(\G)_*$ defined to be $\left(\int_X\ \cdot\  \mu(dx)\right)\circ E_{\G}$. Therefore, it suffices to consider $E_A^M$-groupoids and their $\varphi$-costs to re-prove Gaboriau's results in operator algebra framework with generalizing it to the (even not necessary non-principal) groupoid setting, and indeed many of results in \cite{Gaboriau:InventMath00} can be proved purely in the framework. For example, we can show the following additivity formula of costs of $E_A^M$-groupoids:  

\begin{theorem}\label{Thm3-1} {\rm (cf.~\cite[Th\'{e}o\`{e}me IV.15]{Gaboriau:InventMath00})} Assume that $M$ has a faithful tracial state $\tau \in M_*$ with $\tau\circ E_A^M = \tau$. Let $\mathcal{G}_1 \supseteq \mathcal{G}_3 \subseteq \mathcal{G}_2$ be $E_A^M$-groupoids. Set $N_1 := \mathcal{G}_1''$, $N_2 := \mathcal{G}_2''$ and $N_3 := \mathcal{G}_3''$ {\rm (}all of which clearly contains $A${\rm )}, and let $E_{N_3}^M : M \rightarrow N_3$ be the $\tau$-conditional expectation {\rm (}hence $E_A^M\circ E_{N_3}^M = E_A^M${\rm )}. Suppose that 
\begin{equation*}    
\big(M,E_{N_3}^M\big) = \big(N_1, E_{N_3}^M\big|_{N_1}\big)\underset{N_3}{\bigstar}\big(N_2, E_{N_3}^M\big|_{N_2}\big), 
\end{equation*}
or equivalently $\mathcal{G}_1$, $\mathcal{G}_2$ are $*$-free with amalgamation $N_3$ with respect to $E_{N_3}^M$, and further that $A$ is a MASA in $N_3$ so that $\mathcal{G}_3 = \mathcal{G}\left(N_3\supseteq A\right)$ holds automatically, see the discussion just below Lemma \ref{Lem2-2.1}. {\rm (}Remark here that $A$ needs not to be a MASA in $N_1$ nor $N_2$.{\rm )} Then, if $N_3$ is hyperfinite, then the smallest $E_A^M$-groupoid $\mathcal{G} = \mathcal{G}_1\vee\mathcal{G}_2$ that contains $\mathcal{G}_1$ and $\mathcal{G}_2$ satisfies that  
\begin{equation*} 
C_{\tau}\left(\mathcal{G}\right) = C_{\tau}\left(\mathcal{G}_1\right) + C_{\tau}\left(\mathcal{G}_2\right) - C_{\tau}\left(\mathcal{G}_3\right) 
\end{equation*} 
as long as when $C_{\tau}\left(\mathcal{G}_1\right)$ and $C_{\tau}\left(\mathcal{G}_2\right)$ are both finite. 
\end{theorem} 
This can be regarded as a slight generalization of one of the main results in \cite{Gaboriau:InventMath00} to the groupoid setting. In fact, let $\G$ be a countable discrete measurable groupoid with an invariant probability measure $\mu$, and assume that it is generated by two countable discrete measurable subgroupoids $\G_1$, $\G_2$. If no alternating word in $\G_1\setminus\G_3, \G_2\setminus\G_3$ with $\G_3 := \G_1\cap\G_2$ intersects with the unit space of strictly positive measure, i.e., $\G$ is the ``free product with amalgamation $\G_1 \bigstar_{\G_3} \G_2$" (modulo null set), and $\G_3$ is principal and hyperfinite, then the above formula immediately implies the formula $C_{\mu}(\G) = C_{\mu}\left(\G_1\right)+C_{\mu}\left(\G_2\right)-C_{\mu}\left(\G_3\right)$ as long as when $C_{\mu}\left(\G_1\right)$ and $C_{\mu}\left(\G_2\right)$ are both finite. Here, we need the same task as in \cite{Kosaki:JFA04}. 

Proving the above theorem needs several lemmas and propositions, many of which can be proved based on the essentially same ideas as in \cite{Gaboriau:InventMath00} even in operator algebra framework so that some of their details will be just sketched. 

The next simple fact is probably known but we could not find a suitable reference.  

\begin{lemma}\label{Lem3-2} Let $\mathcal{G}$ be an $E_A^M$-groupoid with $M = \mathcal{G}''$, and assume that $M$ is finite. Then, if $e, f \in A^p$ are equivalent in $M$, denoted by $e \sim_M f$, in the sense of Murray-von Neumann {\rm (}i.e., $l(u) = e$ and $f=r(u)$ for some $u \in M^{pi}${\rm )}, then there is an element $u \in \mathcal{G}$ such that $l(u)=e$ and $r(u)=f$. Hence, under the same assumption, if $p \in A^p$ has the central support projection $c_M(p)=1$, then one can find $v_k \in \mathcal{G}$ in such a way that $\sum_k v_k p v_k^* = 1$.  
\end{lemma} 
\begin{proof} 
The latter assertion clearly follows from the former. Since the linear span of $\mathcal{G}$ becomes a dense $*$-subalgebra of $M$, $e\sim_M f$ implies $eMf\neq\{0\}$ so that there is a $v\in\mathcal{G}$ with $evf\neq0$. Letting $u_0 := fve$ one has $l(u_0) \leq f$ and $r(u_0)\leq e$, and thus $f-l(u_0) \sim_M e-r(u_0)$ since $M$ is finite. Hence, standard exhaustion argument completes the proof. 
\end{proof}    

To prove the next proposition, Gaboriau's original argument still essentially works purely in operator algebra framework. 

\begin{proposition}\label{Prop3-3}{\rm (\cite[Proposition I.9; Proposition I.11]{Gaboriau:InventMath00})} Suppose that $A$ is a MASA in $M$. Then, the following assertions hold true{\rm:} 
\begin{itemize} 
\item[{\rm (a)}] Let $\varphi \in M_*$ be a faithful state with $\varphi\circ E_A^M = \varphi$. If a graphing $\mathcal{U}$ of $\mathcal{G}(M\supseteq A)$ satisfies $C_{\varphi}(\mathcal{U})=C_{\varphi}\left(\mathcal{G}(M\supseteq A)\right)<+\infty$, then $\mathcal{U}$ must be a treeing.
\item[{\rm (b)}] If $M$ is of finite type I  {\rm (}hence $A$ is automatically a Cartan subalgebra{\rm )} and $\tau \in M_*$ is a faithful tracial state {\rm (}n.b., $\tau\circ E_A^M=\tau$ holds automatically{\rm )}, then every treeing $\mathcal{U}$ of $\mathcal{G}(M\supseteq A)$ satisfies that 
\begin{equation*} 
C_{\tau}(\mathcal{U}) = 1 - \tau(e) = C_{\tau}\left(\mathcal{G}(M\supseteq A)\right), 
\end{equation*} 
where $e \in A^p$ is arbitrary, maximal, abelian projection of $M$ {\rm (}hence the central support projection $c_M(e)=1${\rm )}.  
\end{itemize}  
\end{proposition}
\begin{proof} (Sketch) (a) Suppose that $\mathcal{U}$ is not a treeing. Then, one can choose a word $v_{\ell}^{\varepsilon_{\ell}}\cdots v_1^{\varepsilon_1}$ in $\mathcal{U}\sqcup\mathcal{U}^*$ of reduced form in the formal sense in such a way that $E_A^M\left(v_{\ell}^{\varepsilon_{\ell}}\cdots v_1^{\varepsilon_1}\right)\neq0$ but every proper subword $v_i^{\varepsilon_i}\cdots v_j^{\varepsilon_j}$ satisfies that $E_A^M\left(v_i^{\varepsilon_i}\cdots v_j^{\varepsilon_j}\right)=0$. It is plain to find mutually orthogonal nonzero $e_1,\dots,e_{\ell} \in A^p$ with $e_k \leq r\left(v_k^{\varepsilon_k}\right)$ satisfying that $v_k^{\varepsilon_k}e_k v_k^{\varepsilon_k}{}^* = e_{k+1}$ ($k=1,\dots,\ell-1$) and $v_{\ell}^{\varepsilon_{\ell}} e_{\ell} v_{\ell}^{\varepsilon_{\ell}}{}^* = e_1$, where the following simple fact is needed: If $A$ is a MASA in $M$, then any $u \in \mathcal{G}(M\supseteq A)\setminus A^{pi}$ has a nonzero $e \in A^p$ such that $e\leq r(v)$ and $e(vev^*)=0$. Thus, $\mathcal{V} := \mathcal{U}\setminus\left\{v_{\ell}\right\}\sqcup\left\{\left(l\left(v_{\ell}^{\varepsilon_{\ell}}\right)-e_1\right)v_{\ell}^{\varepsilon_{\ell}}\right\}$ becomes a graphing and satisfies $C_{\varphi}(\mathcal{U})\gneqq C_{\varphi}(\mathcal{V})$, a contradiction.  
   
(b) Assume that $M = M_n\left(\mathbf{C}\right)$. Let $\mathcal{V}$ be a graphing of $\mathcal{G}(M\supseteq A)$. Let $p_1,\dots,p_n \in A^p$ be the mutually orthogonal minimal projections in $M$, and define the new graphing $\mathcal{V}'$ to be the collection of all nonzero $p_i v p_j$ with $i,j=1,\dots,n$ and $v \in \mathcal{V}$, each of which is nothing but a standard matrix unit (modulo scalar multiple). Note that $C_{\tau}(\mathcal{V}) = C_{\tau}\left(\mathcal{V}'\right)$ by the construction, and it is plain to see that if $\mathcal{V}$ is a treeing then so is $\mathcal{V}'$ too. We then construct a (non-oriented, geometric) graph whose vertices are $p_1,\dots,p_n$ and whose edges given by $\mathcal{V}'$ with regarding each $p_i v p_j \in \mathcal{V}'$ as an arrow connecting between $p_i$ and $p_j$. It is plain to see that a sub-collection $\mathcal{U}$ of $\mathcal{V}'$ is a treeing of $\mathcal{G}(M\supseteq A)$ if and only if the subgraph whose edges are given by only $\mathcal{U}$ forms a maximal tree. Therefore, a standard fact in graph theory (see e.g.~\cite[\S\S2.3]{Serre:book})  tells that $\mathcal{V}'$ contains a treeing $\mathcal{U}$ of $\mathcal{G}(M\supseteq A)$ or $\mathcal{V}'$ becomes a treeing when so is $\mathcal{V}$ itself. Such a treeing is determined as a collection of matrix units $e_{i_1 j_1},\dots,e_{i_{n-1} j_{n-1}}$ up to scalar multiples with the property that each of $1,\dots,n$ appears at least once in the subindices $i_1,j_1,\dots,i_{n-1},j_{n-1}$. Hence $C_{\varphi}(\mathcal{V})=C_{\tau}\left(\mathcal{V}'\right)\geq C_{\varphi}(\mathcal{U}) = 1-1/n$, which implies the desired assertion in the special case of $M=M_n(\mathbf{C})$. The simultaneous central decomposition of $M \supseteq A$ reduces the general case to the above simplest case we have already treated. Proving that any treeing attains $C_{\tau}\left(\mathcal{G}(M\supseteq A)\right)$ needs the following simple fact: Let $\mathcal{U}$ be a graphing of $\mathcal{G}(M\supseteq A)$, and set $\mathcal{U}(\omega) := \left\{u(\omega) : u \in \mathcal{U}\right\}$ with $u = \int^{\oplus}_{\Omega} u(\omega) d\omega$ in the central decomposition of $M$ with $\mathcal{Z}(M) = L^{\infty}(\Omega) \subseteq A$. Then, $\mathcal{U}$ is $*$-free with respect to $E_A^M$ (or other words, say a treeing) if and only if so is $\mathcal{U}(\omega)$ with respect to $E_{A(\omega)}^{M(\omega)}$ for a.e.~$\omega\in\Omega$ with $E_A^M = \int_{\Omega}^{\oplus} E_{A(\omega)}^{M(\omega)} d\omega$, see e.g.~the proof of \cite[Theorem 5.1]{Ueda:ASPM04}. 
\end{proof} 

\begin{remark}\label{Rem3-4}{\rm  (1) In the above (a), it cannot be avoided to assume that $A$ is a MASA in $M$, that is, the assertion no longer holds true in the non-principal groupoid case. In fact, let $M := L\left(\mathbb{Z}_N\right)$ be the group von Neumann algebra associated with cyclic group $\mathbb{Z}_N$ and $\tau_{\mathbb{Z}_N}$ be the canonical tracial state. Then, $\mathcal{G}(\mathbb{Z}_N) := \mathbb{T}1\cdot\lambda(\mathbb{Z}_N)$ is a $\tau_{\mathbb{Z}_N}(\ \cdot\ )1$-groupoid, and it is trivial that $C_{\tau_{\mathbb{Z}_N}}\left(\mathcal{G}(\mathbb{Z}_N)\right) = C_{\tau_{\mathbb{Z}_N}}\left(\{\lambda(\bar{1})\}\right)$ with the canonical generator $\bar{1} \in \mathbb{Z}_N$. This clearly provides a counter-example. 

(2) Notice that the cost $C_{\tau_G}(\mathcal{G}(G))$ of a group $G$ is clearly the smallest number $n(G)$ of generators of $G$, and hence Theorem \ref{Thm3-1} provides a quite natural formula, that is, $n(G\bigstar H) = n(G)+n(H)$. One should here note that the $\ell^2$-Betti numbers of discrete groupoids (\cite{Gaboriau:IHES02}, and also \cite{Sauer:Preprint03}) recover the group $\ell^2$-Betti numbers when a given groupoid is a group (see e.g.~the approach in \cite{Sauer:Preprint03}). 

(3) Assume that $M$ is properly inifinite and $A$ is a Cartan subalgebra in $M$. Based on the fact that the inclusion $B\big(\ell^2(\mathbb{N)}\big) \supseteq \ell^{\infty}(\mathbb{N})$ can be embedded into $M\supseteq A$, it is not difficult to see that $C_{\varphi}\left(\mathcal{G}(M\supseteq A)\right)=\frac{1}{2}$ for every faithful state $\varphi\in M_*$ with $\varphi\circ E_A^M=\varphi$. Therefore, the idea of costs seems to fit for nothing in the infinite case with general states.  

(4) One of the key ingredients in the proof of (b) can be illustrated by 
\begin{equation*} 
M_3(\mathbf{C}) \cong \begin{bmatrix} * & * & \\ * & * & \\ & & * \end{bmatrix} \underset{\begin{bmatrix} * & & \\ & * & \\ & & * \end{bmatrix}}{\bigstar} \begin{bmatrix} * & & \\ & * & * \\ & * & * \end{bmatrix}   
\end{equation*}
which provides the treeing $e_{12}, e_{23}$ of  $\mathcal{G}(M\supseteq A)$ with $M=M_3(\mathbf{C})$. This kind of facts are probably known, and specialists in free probability theory are much familiar with similar phenomena in the context of (operator) matrix models of semicircular systems.   
}
\end{remark}       

Throughout the rest of this section, let us assume that $\mathcal{G}$ is an $E_A^M$-groupoid with $M=\mathcal{G}''$ and $\tau \in M_*$ is a faithful tracial state with $\tau\circ E_A^M = \tau$. For a given $p \in A^p$ we denote by $p\mathcal{G}p$ the set of $pup$ with $u\in\mathcal{G}$, which becomes an $E_{Ap}^{pMp}$-groupoid with $E_{Ap}^{pMp} := E_A^M\big|_{pMp}$. The next lemma is technical but quite important, and shown in the same way as in Gaboriau's. It is a graphing counterpart of the well-known construction of induced transformations (see e.g.~\cite[p.13--14]{Friedman:Book70}).  

\begin{lemma}\label{Lem3-5}{\rm (cf.~\cite[Lemme II.8]{Gaboriau:InventMath00})} Let $p \in A^p$ be such that the central support projection $c_M(p)=1$, and $\mathcal{U}$ be a graphing of $\mathcal{G}$. Then, there are a treeing $\mathcal{U}_v$ and a graphing $\mathcal{U}_h$ of $p\mathcal{G}p$ with the following properties{\rm :}  
\begin{itemize}
\item[{\rm (a)}] $p$ is an abelian projection of $M_v:=A\vee\mathcal{U}_v{}''$ with $c_{M_v}(p) = 1${\rm ;} 
\item[{\rm (b)}] For a graphing $\mathcal{V}$ of $p\mathcal{G}p$, $\mathcal{U}_v\sqcup\mathcal{V}$ becomes a graphing of $\mathcal{G}${\rm ;} 
\item[{\rm (c)}] For a graphing $\mathcal{V}$ of $p\mathcal{G}p$, $\mathcal{U}_v\sqcup\mathcal{V}$ is a treeing of $\mathcal{G}$ if and only if so is $\mathcal{V}${\rm ;} 
\item[{\rm (d)}] $C_{\tau}(\mathcal{U}) = C_{\tau}(\mathcal{U}_v)+C_{\tau}(\mathcal{U}_h)$ and $C_{\tau}(\mathcal{U}_v) = 1-\tau(p)$. 
\end{itemize}
\end{lemma}     
\begin{proof} (Sketch) Let $\mathcal{U}^{(\ell)}$ be the set of words in $\mathcal{U}\sqcup\mathcal{U}^*$ of reduced form in the formal sense and of length $\ell\geq 1$, and set $q_{\ell} := \bigvee_{w \in \mathcal{U}^{(\ell)}} w p w^*$. Since $A$ is commutative, we can construct inductively the projections $p_{\ell} \in A^p$ by $p_{\ell} := q_{\ell}(1-p_1-\cdots-p_{\ell-1})$ with $p_0:=p$. Letting $p_0 := p$ we have $\sum_{\ell\gneq0} p_{\ell} = 1$ thanks to $c_M(p)=1$. For each $u \in \mathcal{U}$, we define $u_{\ell_1\ell_2} := p_{\ell_1} u p_{\ell_2} \in \mathcal{G}$ with $\ell_1,\ell_2 \in \mathbb{N}\sqcup\{0\}$, and consider the new collection $\widetilde{\mathcal{U}} := \bigsqcup_{\ell_1,\ell_2\geq0}\widetilde{\mathcal{U}}_{\ell_1,\ell_2}$ with $\widetilde{\mathcal{U}}_{\ell_1,\ell_2} := \left\{ u_{\ell_1,\ell_2} : u \in \mathcal{U} \right\}$ instead of the original $\mathcal{U}$ (without changing the $\tau$-costs). Replacing $u_{\ell_1\ell_2}$ by its adjoint if $\ell_2 \lneqq \ell_1$ we may and do assume that $\widetilde{\mathcal{U}}_{\ell_1,\ell_2} = \emptyset$ as long as when $\ell_2 \lneqq \ell_1$. Then, it is not so hard to see that $p_{\ell} = \bigvee_{v\in\widetilde{\mathcal{U}}_{\ell-1,\ell}} r(v)$ for every $\ell \geq 1$. Numbering $\widetilde{\mathcal{U}}_{\ell-1,\ell} = \left\{v_1,v_2,\dots\right\}$ we construct a partition $p_{\ell} = \sum_k s_k$ in $A^p$ inductively by $s_k := r(v_k)(1-s_1-\cdots-s_{k-1})$, and set $\widetilde{\mathcal{U}}_{\ell-1,\ell}' := \left\{v_k s_k\right\}_k$ and $\widetilde{\mathcal{U}}_{\ell-1,\ell}'' := \left\{v_k (1-s_k)\right\}_k$. Set $\mathcal{U}_v := \bigsqcup_{\ell\geq1}\widetilde{\mathcal{U}}'_{\ell-1.\ell}$, and then it is clear that the (right support) projections  $r(v)$, $v\in\mathcal{U}_v$, are mutually orthogonal and moreover that $\sum_{v\in\widetilde{\mathcal{U}}_{\ell-1,\ell}'} r(v) = p_{\ell}$ (hence $\sum_{v\in\mathcal{U}_v}r(v)=1-p$). Set $\mathcal{U}_v^{[k,\ell]} := \big\{ v_{k k+1}\cdots v_{\ell-1 \ell} \neq 0 : v_{j-1 j} \in \widetilde{\mathcal{U}}_{j-1,j}'\big\}$ with $k\lneqq\ell$, and define $\mathcal{U}_h$ to be the collection of elements in $\mathcal{G}$ of the form, either $v \in \mathcal{U}_{0,0}$ or $w_1 v w_2^* \neq 0$ with either $w_1 \in \mathcal{U}_v^{[0,\ell_1]}$, $v \in\widetilde{\mathcal{U}}_{\ell_1,\ell_2}$, $w_2 \in \mathcal{U}_v^{[0,\ell_2]}$ ($\ell_1 = \ell_2$ or  $\ell_1 \leq \ell_2-2$); or $w_1 \in \mathcal{U}_v^{[0,\ell]}$, $v \in\widetilde{\mathcal{U}}_{\ell-1,\ell}''$, $w_2 \in \mathcal{U}_v^{[0,\ell-1]}$. 
It is not so hard to verify that all the assertions (a)-(d) hold for the collections $\mathcal{U}_v$, $\mathcal{U}_h$ that we just constructed. (Note here that the trace property of $\tau$ is needed only for verifying the assertion (d).)         
\end{proof}  
 
\begin{remark}\label{Rem3-5.1} {\rm We should remark that $M_v$ is constructed so that $A$ is a Cartan subalgebra in $M_v$. Let $\mathcal{G}_v$ be the smallest $E_A^M$-groupoid that contains $\mathcal{U}_v$, and hence $N_v = \mathcal{G}_v''$ is clear. By the construction of $\mathcal{U}_v$ one easily see that any non-zero word in $\mathcal{U}_v\sqcup\mathcal{U}_v^*$ must be in either $\mathcal{U}_v^{[k,\ell]}$ or its adjoint set so that $p\mathcal{G}_v p = A^{pi}p$ by Lemma \ref{Lem2-2}. (This pattern of argument is used to confirm that $\mathcal{U}_v$ is a treeing.) Hence, we get $\mathcal{Z}(M_v)p = pM_v p = Ap$, by which with $c_{M_v}(p)=1$ it immediately follows that $A'\cap M_v= A$, thanks to Lemma \ref{Lem3-2}.}  
\end{remark}  

\begin{proposition}\label{Prop3-6}{\rm (cf.~\cite[Proposition II.6]{Gaboriau:InventMath00})} Let $p \in A^p$ be such that the central support projection $c_M(p)=1$. Then, the following hold true{\rm :} 
\begin{itemize}
\item $\mathcal{G}$ is treeable if and only if so is $p\mathcal{G}p${\rm ;} 
\item $C_{\tau}(\mathcal{G}) - 1 = C_{\tau|_{pMp}}(p\mathcal{G}p)-\tau(p)$. 
\end{itemize}
\end{proposition}  
\begin{proof} The first assertion is nothing less than Lemma \ref{Lem3-5} (c). The second is shown as follows. By Lemma \ref{Lem3-5} (d), we have $C_{\tau}(\mathcal{U}) \geq C_{\tau}\left(p\mathcal{G}p\right) + 1-\tau(p)$ for every graphing $\mathcal{U}$ of $\mathcal{G}$ so that $C_{\tau}(\mathcal{G}) -1 \geq C_{\tau}\left(p\mathcal{G}p\right)-\tau(p)$. Let $\varepsilon>0$ be arbitrary small. Choose a graphing $\mathcal{V}_{\varepsilon}$ so that $C_{\tau}\left(\mathcal{V}_{\varepsilon}\right)\leq C_{\tau}\left(p\mathcal{G}p\right)+\varepsilon$. With $\mathcal{U}_v$ as in Lemma \ref{Lem3-5} the new collection $\mathcal{U}_{\varepsilon} :=\mathcal{U}_v\sqcup\mathcal{V}_{\varepsilon}$ becomes a graphing of $\mathcal{G}$ by Lemma \ref{Lem3-5} (b), and hence $C_{\tau}(\mathcal{G})\leq C_{\tau}\left(\mathcal{U}_{\varepsilon}\right)=1-\tau(p)+C_{\tau}\left(\mathcal{V}_{\varepsilon}\right)$ by Lemma \ref{Lem3-5} (d). Hence, $C_{\tau}(\mathcal{G})-1\leq C_{\tau}\left(\mathcal{V}_{\varepsilon}\right)-\tau(p)\leq C_{\tau}\left(p\mathcal{G}p\right)+\varepsilon-\tau(p)\searrow C_{\tau}\left(p\mathcal{G}p\right)-\tau(p)$ as $\varepsilon\searrow0$. 
\end{proof}   

\begin{corollary}\label{Cor3-7} {\rm (\cite[Proposition 1, Theorem 2]{Levitt:ErgodicTheoryDynamSystems95},\cite[Proposition III.3, Lemme III.5]{Gaboriau:InventMath00})} 
{\rm (a)} Assume that $M$ is of type II$_1$ and $A$ is a Cartan subalgebra in $M$. Then, $C_{\tau}\left(\mathcal{G}(M\supseteq A)\right)\geq1$, and the equality holds if $M$ is further assumed to be hyperfinite. 

{\rm (b)} Assume that $M$ is hyperfinite and $A$ is a Cartan subalgebra in $M$. Then, every treeing $\mathcal{U}$ of $\mathcal{G}(M\supseteq A)$ {\rm (}it always exists{\rm )} satisfies that 
\begin{equation*} 
C_{\tau}(\mathcal{U}) = 1-\tau(e) = C_{\tau}\left(\mathcal{G}(M\supseteq A)\right),
\end{equation*} 
where $e \in A^p$ is arbitrary, maximal abelian projection of $M$ {\rm (}hence the central support projection $c_M(e)$ must coincide with that of type I direct summand{\rm )}. 

{\rm (c)} Let $N$ be a hyperfinite intermediate von Neumann subalgebra between $M\supseteq A$, and assume that $A$ is a Cartan aubalgebra in $N$. Let $\mathcal{U}$ be a treeing of $\mathcal{G}(N\supseteq A)$ and suppose that $\mathcal{G}$ contains $\mathcal{G}(N\supseteq A)$. Then, for each $\varepsilon>0$, there is a graphing $\mathcal{U}_{\varepsilon}$ of $\mathcal{G}$ enlarging $\mathcal{U}$ such that $C_{\tau}\left(\mathcal{U}_{\varepsilon}\right) \leq C_{\tau}(\mathcal{G}) + \varepsilon$.  
\end{corollary}
\begin{proof} (a) It is known that for each $n \in \mathbb{N}$ there is an $n\times n$ matrix unit system $e_{ij} \in \mathcal{G}(M\supseteq A)$ ($i,j=1,\dots,n$) such that all $e_{ii}$'s are chosen from $A^p$. Then, Proposition \ref{Prop3-6} implies that $C_{\tau}(\mathcal{G})=C_{\tau|_{e_{11}Me_{11}}}\left(e_{11}\mathcal{G}e_{11}\right)+1-\tau\left(e_{11}\right)\geq1-\tau\left(e_{11}\right)=1-1/n\nearrow1$ as $n\rightarrow\infty$. The equality in the hyperfinite case clearly follows from celebrated Connes, Feldman and Weiss' theorem \cite{ConnesFeldmanWeiss:ErgodicTh81} (also \cite{Popa:MathScand85} for its operator algebraic proof).  

(b) Choose an incereasing sequence of type I von Neumann subalgebras $A \subseteq M_1 \subseteq \cdots \subseteq M_k \nearrow M$. By Dye's lemma (or Lemma \ref{Lem2-2.1}), each $u \in \mathcal{U}$ has a unique projection $e_k(u) \in A^p$ such that $e_k(u)\leq l(u)$ and $E_{N_k}^M(u)=e_k(u)u$, where $E_{M_k}^M : M \rightarrow M_k$ is the $\tau$-conditional expectation. Set $\mathcal{U}_k := \left\{ e_k(u)u  : u \in \mathcal{U}\right\}$ and $N_k := A\vee\mathcal{U}_k''$ being of type I. Clearly, each $\mathcal{U}_k$ is a treeing of $\mathcal{G}(N_k\supseteq A)$, and hence Proposition \ref{Prop3-3} (b) says that  $C_{\tau}\left(\mathcal{U}_k\right)=C_{\tau}\left(\mathcal{G}\left(N_k\supseteq A\right)\right)=1-\tau\left(e_k\right)$ for every maximal abelian projection $e_k \in A^p$ of $N_k$. The non-commutative Martingale convergence theorem (e.g.~\cite[Lemma 2]{Connes:JFA75}) shows that $e_k(u) u = E_{M_k}^M(u) \rightarrow u$ in $\sigma$-strong* topology, as $k\rightarrow\infty$, for every $u \in \mathcal{U}$. Hence, we get $C_{\tau}\left(\mathcal{U}\right)=\lim_{k\rightarrow\infty}C_{\tau}\left(\mathcal{U}_k\right)=\lim_{k\rightarrow\infty}C_{\tau}\left(\mathcal{G}\left(N_k\supseteq A\right)\right)$ and $N_k = A\vee\mathcal{U}_k''  \nearrow A\vee\mathcal{U}''=M$. Let $e \in A^p$ be a maximal abelian projection of $M$. Then, Proposition \ref{Prop3-3} (b) and the above (a) show that $C_{\tau}\left(\mathcal{G}(M\supseteq A)\right)=1-\tau(e)$. Since $e$ must be an abelian projection of each $N_k$, one can choose $e_1,e_2,\dots \in A^p$ in such a way that each $e_k$ is a maximal abelian projection of $N_k$ and greater than $e$. It is standard to see that $e=\bigwedge_{k=1}^{\infty} e_k$ so that $\tau(e_k) \geq \tau\big(\bigwedge_{k'=1}^k e_{k'}\big) \searrow \tau(e)$ as $k\rightarrow\infty$. Therefore, $C_{\tau}(\mathcal{U})=\lim_{k\rightarrow\infty}C_{\tau}\left(\mathcal{G}\big(N_k\supseteq A\right)\big)=\lim_{k\rightarrow\infty}\left(1-\tau\left(e_k\right)\right)\leq\lim_{k\rightarrow\infty}\left(1-\tau\left(\bigwedge_{k'=1}^k e_{k'}\right)\right) = 1-\tau(e)$, and then it follows immediately that $C_{\tau}(\mathcal{U}) = 1-\tau(e) = C_{\tau}\left(\mathcal{G}(M\supseteq A)\right)$.    

(c) Let $N = N_{\mathrm{I}}\oplus N_{\mathrm{II}_1} \supseteq A = A_{\mathrm{I}}\oplus A_{\mathrm{II}_1}$ be the decomposition into the finite type I and the type II$_1$ parts. Looking at the decomposition, one can find a projection $p_{\varepsilon}=p_I\oplus p_{\mathrm{II}_1}^{\varepsilon} \in A^p$ in such a way that $p_{\mathrm{I}}$ is an abelian projection of $N_{\mathrm{I}}$ with $c_{N_{\mathrm{I}}}(p_{\mathrm{I}}) = 1_{N_{\mathrm{I}}}$ and $\tau\left(p_{\mathrm{II}_1}^{\varepsilon}\right)<\varepsilon/2$ with $c_{N_{\mathrm{II}_1}}\left(p_{\mathrm{II}_1}^{\varepsilon}\right)=1_{N_{\mathrm{II}_1}}$. Choose a graphing $\mathcal{V}_{\varepsilon}$ of $p_{\varepsilon}\mathcal{G}p_{\varepsilon}$ in such a way that $C_{\tau}\left(\mathcal{V}_{\varepsilon}\right)\leq C_{\tau}\left(p_{\varepsilon}\mathcal{G}p_{\varepsilon}\right)+\varepsilon/2$, and then set $\mathcal{U}_{\varepsilon} := \mathcal{U}\sqcup\mathcal{V}_{\varepsilon}$. Since $c_{N}\left(p_{\varepsilon}\right)=1$, $\mathcal{U}_{\varepsilon}$ is a graphing $\mathcal{G}$ thanks to Lemma \ref{Lem3-2}. Then, Lemma \ref{Lem3-5} (b) implies that $C_{\tau}\left(\mathcal{V}_{\varepsilon}\right) 
\leq C_{\tau}\left(p_{\varepsilon}\mathcal{G}p_{\varepsilon}\right)+\varepsilon/2 = 
C_{\tau}(\mathcal{G})-1+\tau\left(p_{\varepsilon}\right)+\varepsilon/2 = C_{\tau}(\mathcal{G})-\left(1-\tau\left(p_{\mathrm{I}}\right)\right)+\tau\left(p_{\mathrm{II}_1}^{\varepsilon}\right)+\varepsilon/2$, 
and thus by Proposition \ref{Prop3-3} (b) we get $C_{\tau}\left(\mathcal{V}_{\varepsilon}\right) \leq C_{\tau}(\mathcal{G}) - C_{\tau}(\mathcal{U}) + \varepsilon$, which implies the desired assertion.      
\end{proof} 

\begin{remark}\label{Rem3-8} {\rm (1) The proof of (b) in the above also shows ``hyperfinite monotonicity," which asserts as follows. Assume that $M$ is hyperfinite and $A$ is a Cartan subalgebra in $M$. For any intermediate von Neumann subalgebra $N$ between $M \supseteq A$ (in which $A$ becomes automatically a Cartan subalgebra thanks to Dye's lemma, see the discussion above Lemma \ref{Lem2-3}), we have $C_{\tau}\left(\mathcal{G}(N\supseteq A)\right) \leq C_{\tau}\left(\mathcal{G}(M\supseteq A)\right)$. Furthermore, we have $\lim_{k\rightarrow\infty}C_{\tau}\big(\mathcal{G}(M_k\supseteq A)\big)=C_{\tau}\big(\mathcal{G}(M\supseteq A)\big)$ for any increasing sequence $A \subseteq M_1 \subseteq M_2 \subseteq \cdots \subseteq M_k \nearrow M$ of von Neumann subalgebras. Note that this kind of fact on free entropy dimension was provided by K.~Jung \cite{Jung:TAMS03}. 

(2) Related to (c) one can show the following (cf.~\cite[Lemme V.3]{Gaboriau:InventMath00}): Let $u \in \mathcal{G}$, and $\mathcal{G}_0 \subseteq \mathcal{G}$ be an $E_A^M$-groupoid, and set $N := \big( r(u)\mathcal{G}_0r(u)\big)''\vee\big(u^*\mathcal{G}_0u\big)''$. If $e \in \big(Ar(u)\big)^p$ has $c_N(e) = r(u)$, then $\mathcal{G}_0\vee\{u\} = \mathcal{G}_0\vee\{ue\}$ so that $C_{\tau}(\mathcal{G}_0\vee\{u\}) \leq C_{\tau}(\mathcal{G}_0) + \tau(e)$. Here, ``$\vee$" means the symbol of generation as $E_A^M$-groupoid. In fact, by Lemma \ref{Lem3-2} one finds $v_k \in r(u)\mathcal{G}_0 r(u) \vee u^*\mathcal{G}_0 u$ so that $\sum_k v_k e v_k^* = r(u)$. Since $v_k \in u^*\mathcal{G}_0 u$, one has $v_k = u^* w_k u$ for some $w_k \in l(u)\mathcal{G}_0 l(u)$ so that $\sum_k w_k(ue) v_k^* = u$.  This fact can be used in many actual computations, and in fact it tells us that the cost of an $E_A^M$-groupoid can be estimated by that of its ``normal $E_A^M$-subgroupoid" with a certain condition. Its free entropy dimension counterpart seems an interesting question. }
\end{remark}  

\begin{proposition}\label{Prop3-9} Let $\mathcal{G}_1\supseteq\mathcal{G}_3\subseteq\mathcal{G}_2$ be $E_A^M$-groupoids, and let $\mathcal{G} = \mathcal{G}_1\vee\mathcal{G}_2$ be the smallest $E_A^M$-groupoid that contains $\mathcal{G}_1$, $\mathcal{G}_2$. If $\mathcal{G}_3''$ is hyperfinite and if $A$ is a MASA in $\mathcal{G}_3''$ {\rm (}and hence $\mathcal{G}_3=\mathcal{G}\left(\mathcal{G}_3''\supseteq A\right)$ is automatic{\rm )}, then 
\begin{equation*} 
C_{\tau}(\mathcal{G}) \leq C_{\tau}\left(\mathcal{G}_1\right)+C_{\tau}\left(\mathcal{G}_2\right)-C_{\tau}\left(\mathcal{G}_3\right). 
\end{equation*} 
\end{proposition} 
\begin{proof} Choose a treeing $\mathcal{U}$ of $\mathcal{G}_3$ so that $C_{\tau}\left(\mathcal{G}_3\right)=C_{\tau}(\mathcal{U})$ by Corollary \ref{Cor3-7} (b). Let $\varepsilon>0$ be arbitrary small. By Corollary \ref{Cor3-7} (c), one can choose graphings $\mathcal{U}_{\varepsilon}^{(i)}$ of $\mathcal{G}_i$ enlarging $\mathcal{U}$, $i=1,2$, so that $C_{\tau}\big(\mathcal{U}_{\varepsilon}^{(i)}\big)\leq C_{\tau}\left(\mathcal{G}_i\right)+\varepsilon/2$. Thus, $C_{\tau}(\mathcal{G})\leq C_{\tau}\left(\big(\mathcal{U}_{\varepsilon}^{(1)}\setminus\mathcal{U}\big)\sqcup\big(\mathcal{U}_{\varepsilon}^{(2)}\setminus\mathcal{U}\big)\sqcup\mathcal{U}\right) = C_{\tau}\big(\mathcal{U}_{\varepsilon}^{(1)}\big)+C_{\tau}\big(\mathcal{U}_{\varepsilon}^{(2)}\big)-C_{\tau}(\mathcal{U})\leq C_{\tau}\left(\mathcal{G}_1\right)+C_{\tau}\left(\mathcal{G}_2\right)-C_{\tau}\left(\mathcal{G}_3\right)+\varepsilon\searrow C_{\tau}\left(\mathcal{G}_1\right)+C_{\tau}\left(\mathcal{G}_2\right)-C_{\tau}\left(\mathcal{G}_3\right)$ as $\varepsilon\searrow0$. 
\end{proof}

To prove Theorem \ref{Thm3-1}, it suffices to show the inequality $C_{\tau}(\mathcal{G})\geq C_{\tau}\left(\mathcal{G}_1\right)+C_{\tau}\left(\mathcal{G}_2\right)-C_{\tau}\left(\mathcal{G}_3\right)$ thanks to Proposition \ref{Prop3-9}. To do so, we begin by providing a simple fact on general amalgamated free products of von Neumann algebras. 

\begin{lemma}\label{Lem3-10} Let 
\begin{equation*} 
\big(N,E_{N_3}^N\big)=
\big(N_1,E_{N_3}^{N_1}\big)\underset{N_3}{\bigstar}\big(N_2,E_{N_3}^{N_2}\big)
\end{equation*} 
be an amalgamated free product of {\rm (}$\sigma$-finite{\rm )} von Neumann algebras, and $L_1$ and $L_2$ be von Neumann subalgebras of $N_1$ and $N_2$, respectively. Suppose that $\begin{matrix} 
N_i & \supset & L_i \\ 
\cup & & \cup \\
N_3 & \supset & N_3\cap L_i
\end{matrix}$ has faithful normal conditional expectations 
\begin{equation*} 
E_{L_i}^{N_i} : N_i \rightarrow L_i, \quad 
E_{N_3 \cap L_i}^{N_3} : N_3 \rightarrow N_3 \cap L_i, \quad 
E_{N_3 \cap L_i}^{L_i} : L_i \rightarrow N_3 \cap L_i, 
\end{equation*} 
and form commuting squares {\rm (}see e.g.~{\rm \cite[p.~513]{EvansKawahigashi:Book})} for both $i=1,2$. If $L_1\cap N_3 = L_2\cap N_3$ and further $N = L_1\vee L_2$ as von Neumann algebra, then $L_1 = N_2$ and $L_2 = N_2$ must hold true. 
\end{lemma} 
\begin{proof} Note that the amalgamated free product 
\begin{equation*} 
\big(L, E_{L_3}^L\big) = 
\big(L_1, E_{L_3}^L\big)\underset{L_1\cap N_3 = L_2\cap N_3}{\bigstar}\big(L_2, E_{L_3}^{L_2}\big) 
\end{equation*} 
can be naturally embedded into $\big(N, E_{N_3}^N\big)$ thanks to the commuting square assumption. Then, it is plain to see that $\begin{matrix} 
N & \supset & L \\ 
\cup & & \cup \\
N_i & \supset & L_i
\end{matrix}$ 
form commuting squares too, i.e., $E_L^N\big|_{N_i} = E_{L_i}^{N_i}$, $i=1,2$, by which the desired assertion is immediate. 
\end{proof}  

The next technical lemma plays a key r\^{o}le in the proof of Theorem \ref{Thm3-1}. 

\begin{lemma}\label{Lem3-11} {\rm (\cite[IV.37]{Gaboriau:InventMath00})} Assume the same setup as in Theorem \ref{Thm3-1}. Let $\mathcal{V} = \mathcal{V}_1\sqcup\mathcal{V}_2$ be a graphing of $\mathcal{G}$ with collections $\mathcal{V}_1$, $\mathcal{V}_2$ of elements in $\mathcal{G}_1$, $\mathcal{G}_2$, respectively. Then, one can construct two collections $\mathcal{V}'_1$, $\mathcal{V}'_2$ of elements in $\mathcal{G}_3 = \mathcal{G}\left(N_3\supseteq A\right)$ in such a way that 
\begin{itemize} 
\item[{\rm (i)}] $\mathcal{V}'_1\sqcup\mathcal{V}'_2$ is a treeing{\rm ;} 
\item[{\rm (ii)}] $\mathcal{V}_i\sqcup\mathcal{V}'_i$ is a graphing of $\mathcal{G}_i$ for $i=1,2$, respectively. 
\end{itemize} 
\end{lemma} 
   
Before giving the proof, we illustrate the idea in a typical example. Assume that $M = N_1\bigstar_{N_3}N_2 \supseteq A$ is of the form: 
$N_1 := N_1^{(0)}\otimes M_2(\mathbf{C})\otimes M_2(\mathbf{C})$, $N_2 := N_2^{(0)}\otimes M_2(\mathbf{C})\otimes M_2(\mathbf{C})$, $N_3 := N_3^{(0)}\otimes M_2(\mathbf{C})\otimes M_2(\mathbf{C})$ and their common subalgebra $A := A_0\otimes\mathbf{C}^2\otimes\mathbf{C}^2$, 
where $A_0$ is a common Cartan subalgbera of $N_i^{(0)}$, $i=1,2,3$. Denote by $e_{ij}^{(1)}\otimes e_{k\ell}^{(2)}$, $i,j,k,\ell=1,2$, the standard matrix units in $M_2(\mathbf{C})\otimes M_2(\mathbf{C})$. Let $\mathcal{V}_i^{(0)}$ be a graphing of $\mathcal{G}_i := \mathcal{G}\big(N_i^{(0)}\supseteq A_0\big)$, $i=1,2$, and set 
\allowdisplaybreaks{
\begin{align*} 
\mathcal{V}_1 &:= \big\{ v\otimes1\otimes1 : v \in \mathcal{V}_1^{(0)} \big\}\sqcup\big\{1\otimes e_{12}^{(1)}\otimes1 \big\}, \\ 
\mathcal{V}_2 &:= \big\{ v\otimes1\otimes1 : v \in \mathcal{V}_2^{(0)} \big\}\sqcup\big\{1\otimes1\otimes e_{12}^{(1)} \big\}.
\end{align*} 
}Clearly, $\mathcal{V} := \mathcal{V}_1\sqcup\mathcal{V}_2$ becomes a graphing of $\mathcal{G} := \mathcal{G}_1\vee\mathcal{G}_2$. In this example, the collections $\mathcal{V}'_1$, $\mathcal{V}'_2$ in the lemma can be chosen for example to be $\big\{1\otimes e_{11}^{(1)}\otimes e_{12}^{(2)}\big\}$, $\big\{1\otimes e_{12}^{(1)}\otimes1\big\}$, respectively. The proof below goes along the line of this procedure with the help of Lemma \ref{Lem3-10}. 

\begin{proof} 
By Lemma \ref{Lem3-10} with the aid of Lemma \ref{Lem2-2.1} (needed to confirm the required commuting square condition, see the explanation after the lemma) the original (ii) is reduced to showing (ii') $N_3\cap L_1 = N_3\cap L_2$ with $L_i:=A\vee\left(\mathcal{V}_i\sqcup\mathcal{V}'_i\right)''$, $i=1,2$. Choose an increasing sequence of type I von Neumann subalgebras $N_3^{(0)} := A \subseteq N_3^{(1)} \subseteq N_3^{(2)} \subseteq \cdots \subseteq N_3^{(k)} \nearrow N_3$. Let us construct inductively two sequences of collections $\mathcal{V}_1^{(k)}$, $\mathcal{V}_2^{(k)}$ of elements in $\mathcal{G}\big(N_3^{(k)}\supseteq A\big)$ in such a way that 
\begin{itemize} 
\item[(a)] $\mathcal{V}_i^{(k)} \subseteq \mathcal{V}_i^{(k+1)}$ for every $k$ and each $i=1,2$; 
\item[(b)] $\mathcal{V}_1^{(k)}\sqcup\mathcal{V}_2^{(k)}$ is a treeing for every $k$;  
\item[(c)] letting $L_i^{(k)} := A\vee\left(\mathcal{V}_i\sqcup\mathcal{V}_i^{(k)}\right)''$, $i=1,2$, we have 
\begin{itemize}
\item[(1)] $\mathcal{V}_1\sqcup\mathcal{V}_2^{(k)} \subseteq L_1^{(k)}\cap L_2^{(k)}$ for every $k$, 
\item[(2)] $N_3^{(k)}\cap L_1^{(n)} \subseteq L_2^{(k)}$ for every even $k$; 
\item[(3)] $N_3^{(k)}\cap L_2^{(k)} \subseteq L_1^{(k)}$ for every odd $k$. 
\end{itemize} 
\end{itemize} 
((c-2) is needed only for the inductive precedure.) If such collections were constructed, then $\mathcal{V}'_i := \bigcup_k \mathcal{V}_i^{(k)}$, $i=1,2$, would be desired ones. In fact, any word in $\big(\mathcal{V}'_1\sqcup\mathcal{V}'_2\big)\sqcup\big(\mathcal{V}'_1\sqcup\mathcal{V}'_2\big)^*$ of reduced form in the formal sense is in turn one in $\big(\mathcal{V}_1^{(k)}\sqcup\mathcal{V}_2^{(k)}\big)\sqcup\big(\mathcal{V}_1^{(k)}\sqcup\mathcal{V}_2^{(k)}\big)^*$ for some finite $k$ thanks to (a), and thus (i) follows from (b). For each pair $k_1, k_2$, the above (c-2), (c-3) imply, with $k_1, k_2 \leq 2k$, that 
\begin{gather*} 
N_3^{(k_1)}\cap L_1^{(k_2)} \subseteq N_3^{(2k)}\cap L_1^{(2k)} \subseteq N_3^{(2k)}\cap L_2^{(2k)} \subseteq N_3\cap L_2; \\
N_3^{(k_1)}\cap L_2^{(k_2)} \subseteq N_3^{(2k+1)}\cap L_2^{(2k+1)} \subseteq N_3^{(2k+1)}\cap L_1^{(2k+1)} \subseteq N_3\cap L_1. 
\end{gather*}
Hence,  
\begin{equation*} 
\overline{\bigcup_{k_1,k_2} N_3^{(k_1)}\cap L_1^{(k_2)}}^{\text{$\sigma$-s}} \subseteq N_3\cap L_2, \quad \overline{\bigcup_{k_1,k_2} N_3^{(k_1)}\cap L_2^{(k_2)}}^{\text{$\sigma$-s}} \subseteq N_3\cap L_1. 
\end{equation*}
Note here that $N_3\cap L_i = \overline{\bigcup_{k_1,k_2} N_3^{(k_1)}\cap L_i^{(k_2)}}^{\text{$\sigma$-s}}$ for both $i=1,2$, since all 
$\begin{matrix} 
N_i & \supset & L_i \\
\cup & & \cup \\
N_3 & \supset & N_3\cap L_i, 
\end{matrix}\  
\begin{matrix} 
N_i & \supset & L_i \\
\cup & & \cup \\
N_3^{(k)} & \supset & N_3^{(k)}\cap L_i, 
\end{matrix}\ 
\begin{matrix} 
N_i & \supset & L_i^{(k_2)} \\
\cup & & \cup \\
N_3^{(k_1)} & \supset & N_3^{(k_1)}\cap L_i^{(k_2)} 
\end{matrix}$ 
form commuting squares for every $k,k_1,k_2$ and each $i=1,2$, thanks to Dye's lemma (or Lemma \ref{Lem2-2.1}); note here that $A$ is assumed to be a Cartan subalgebra in $N_3$. Then, (ii') follows immediately. 

Set $\mathcal{V}_1^{(0)} = \mathcal{V}_2^{(0)} := \emptyset$. Assume that we have already constructed $\mathcal{V}_1^{(j)}$, $\mathcal{V}_2^{(j)}$, $j=1,2,\dots,k$, and that the next $k+1$ is even (the odd case is also done by the same way). Consider 
\begin{equation*}  
K_1 := N_3^{(k+1)}\cap L_1^{(k)} \supseteq K_0 := N_3^{(k+1)}\cap L_1^{(k)} \cap L_2^{(k)} \left(\supseteq A\right),  
\end{equation*} 
which are clearly of finite type I. Then, one can choose an abelian projection $p \in A^p$ of $K_0$ with the central support projection $c_{K_0}(p)=1$. Also, one can find a treeing $\mathcal{U}_p$ of $\mathcal{G}\left(pK_1 p \supseteq Ap\right)$, see the proof of  Proposition \ref{Prop3-3} (b). Set $\mathcal{V}_1^{(k+1)} := \mathcal{V}_1^{(k)}$, $\mathcal{V}_2^{(k+1)} := \mathcal{V}_2^{(k)}\sqcup\mathcal{U}_p$, which are desired ones in the $k+1$ step. In fact, (a) is trivial, and (b) follows from the (rather trivial) fact that $K_0$ and $pK_1 p$ are $*$-free with respect to $E_A^M$. Note that $N_3^{(k+1)}\cap L_1^{(k+1)} = N_3^{(k+1)}\cap L_1^{(k)} = K_1 = K_0\vee pK_1 p \subseteq L_2^{(k)}\vee\mathcal{U}_p'' = L_2^{(k+1)}$ (by Lemma \ref{Lem3-2} it follows from $c_{K_0}(p)=1$  that $K_1 = K_0\vee pK_1 p$), which is nothing but (c-2). Finally, (c-1) follows from the assumption of induction together with that $\mathcal{U}_p \subseteq pK_1 p \subseteq L_1^{(k)}$.  
\end{proof}

One of the important ideas in Gaboriau's argument is the use of ``adapted systems." It is roughly translated to amplification/reduction procedure in operator algebra framework. 

\begin{proof}{[Proof of Theorem \ref{Thm3-1}]} (Step I: Approximation) By Proposition \ref{Prop3-9}, it suffices to show that $C_{\tau}(\mathcal{G}) \geq C_{\tau}\left(\mathcal{G}_1\right)+C_{\tau}\left(\mathcal{G}_2\right)-C_{\tau}\left(\mathcal{G}_3\right)$ modulo ``arbitrary small error." Let $\varepsilon>0$ be arbitrary small. There is a graphing $\mathcal{V}_{\varepsilon}$ of $\mathcal{G}$ with $C_{\tau}\big(\mathcal{V}_{\varepsilon}\big) \leq C_{\tau}(\mathcal{G}) + \varepsilon/3$, and we choose a graphing $\mathcal{U}_i$ of $\mathcal{G}_i$ with $C_{\tau}\left(\mathcal{U}_i\right)<+\infty$ (thanks to $C_{\tau}\left(\mathcal{G}_i\right) < +\infty$) for each $i=1,2$, and set $\mathcal{U} := \mathcal{U}_1\sqcup\mathcal{U}_2$. Since $C_{\tau}\left(\mathcal{U}\right)=\sum_{u\in\mathcal{U}}\tau\left(l(u)\right)<+\infty$, there is a finite sub-collection $\mathcal{U}_0$ of $\mathcal{U}$ with $\sum_{u\in\mathcal{U}\setminus\mathcal{U}_0} \tau\left(l\left(u\right)\right)\leq\varepsilon/3$. Since both $\mathcal{V}$ and $\mathcal{U}$ are graphings of $\mathcal{G}$, we may and do assume, by cutting each $v \in \mathcal{V}$ by suitable projections in $A^p$ based on Lemma \ref{Lem2-2}, that each $v \in \mathcal{V}$ has a word $w(v)$ in $\mathcal{U}\sqcup\mathcal{U}^*$ of reduced form in the formal sense and a $a(v) \in A^{pi}$ with $v=a(v)w(v)$. Denote by $w_0 := 1, w_1, w_2,\dots$ the all words in $\mathcal{V}\sqcup\mathcal{V}^*$ of reduced form, and by Lemma \ref{Lem2-2} again, each $u\in\mathcal{U}$ is described as $u = \sum_j p_k(u) a_k(u) w_k$ in $\sigma$-strong topology, where $l(u) = \sum_k p_k(u)$ in $A^p$, the $a_k(u)$'s are in $A^{pi}$, and $p_k(u) u = p_k(u) a_k(u) w_k$ for every $k$. Then, we can choose a $k_0\in\mathbb{N}$ (depending only on the finite collection $\mathcal{U}_0$) in such a way that $\sum_{k\geq k_0+1} \tau(p_k(u))\leq\varepsilon/3\sharp(\mathcal{U}_0)$ for all $u\in\mathcal{U}_0$. Set 
\allowdisplaybreaks{
\begin{align*}
\mathcal{X} &:= \left\{ v \in \mathcal{V}_{\varepsilon} : \text{$v$ appears in $w_1,\dots,w_{k_0}$}\right\}; \\
\mathcal{Y} &:= \left\{ p(u) u : u\in\mathcal{U}_0 \right\}\sqcup\left(\mathcal{U}\setminus\mathcal{U}_0\right)
\end{align*} 
}with $p(u) := \sum_{k\geq k_0+1} p_k(u)$ for $u \in \mathcal{U}_0$. Since $u = \sum_{j=0}^{k_0} p_k(u)a_k(u) w_k + p(u)u$ for all $u\in\mathcal{U}_0$, the coollection $\mathcal{Z} := \mathcal{X}\sqcup\mathcal{Y}$ becomes a graphing of $\mathcal{G}$. We have 
\begin{equation}\label{eq1-Thm3-1}
C_{\tau}(\mathcal{Z}) 
= C_{\tau}\left(\mathcal{X}\right)+C_{\tau}\left(\mathcal{Y}\right) 
\leq C_{\tau}\left(\mathcal{V}_{\varepsilon}\right)+C_{\tau}\left(\mathcal{Y}\right)
\leq C_{\tau}(\mathcal{G})+\varepsilon. 
\end{equation}
Clearly, $\mathcal{Y}$ is decomposed into two collections $\mathcal{Y}_1$, $\mathcal{Y}_2$ of elements in $\mathcal{G}_1$, $\mathcal{G}_2$, respectively, as $\mathcal{Y} = \mathcal{Y}_1\sqcup\mathcal{Y}_2$, while $\mathcal{X}$ not in general. Thus, we replace $\mathcal{X}$ by a new ``decomposable" graphing in a sufficiently large amplification of $M\supseteq A$ to use Lemma \ref{Lem3-11}. 

(Step II: Adapted system/Amplification) Notice that each $v \in \mathcal{X} \big(\subseteq \mathcal{V}_{\varepsilon}\big)$ is described as 
\begin{equation*} 
v=a(v)w(v)=a(v)u_{n(v)}(v)^{\delta_{n(v)}(v)}\cdots u_1(v)^{\delta_1(v)}, 
\end{equation*}
where $n(v) \in \mathbb{N}$ and $u_i(v) \in \mathcal{U}$, $\delta_i(v) \in \{1,*\}$ ($i=1,\dots,n(v)$). Cutting each $u_i(v)$ by a suitable projection in $A^p$ and replacing $u_i(v)$ by $u_i(v)^*$ if $\delta_i(v)=*$, etc., we may and do assume the following: $r\left(u_{i+1}(v)\right)=l\left(u_i(v)\right)$ ($i=1,\dots,n(v)-1$); $l(v)=l\left(a(v)\right) \left(=r\left(a(v)\right)\right)=l\left(u_{n(v)}(v)\right)$ and $r(v)=r\left(u_1(v)\right)$; and each $u_i(v)$ is of the form either $eu$ or $eu^*$ with $e \in A^p$, $u \in\mathcal{U}$. In what follows, we ``reveal" all words $u_{n(v)}(v)\cdots u_1(v)$'s as follows. Set $n:=1+\sum_{v\in\mathcal{X}}\left(n(v)-1\right)<+\infty$, and choose a partition $\{2,\dots,n\} = \bigsqcup_{v\in\mathcal{X}}I_v$. Denote by $e_{ij}$ the standard matrix units in $M_n(\mathbf{C})$, by $\mathrm{Tr}_n$ the canonical non-normalized trace on $M_n(\mathbf{C})$, and by $E_{\mathbf{C}^n}^{M_n(\mathbf{C})}$ the $\mathrm{Tr}_n$-conditional expectation from $M_n(\mathbf{C})$ onto the diagonal matrices $\mathbf{C}^n \subseteq M_n(\mathbf{C})$. Let $M^{n} := M\otimes M_n(\mathbf{C}) \supseteq A^{n} := A\otimes\mathbf{C}^n$ be the $n$-amplifications and write $\tau^{n} := \tau\otimes\mathrm{Tr}_n \in M^{n}_*$. For each $v \in \mathcal{X}$, we define the $n(v)$ elements $\tilde{u}_1(v),\dots,\tilde{u}_{n(v)}(v) \in \mathcal{G}\left(M^{n} \supseteq A^{n}\right)$ by 
\allowdisplaybreaks{
\begin{align*} 
\tilde{u}_1(v) &:= u_1(v)\otimes e_{i_2 1}, \\
\tilde{u}_2(v) &:= u_2(v)\otimes e_{i_3 i_2}, \\ 
&\phantom{aaaa}\vdots \\
\tilde{u}_{n(v)}(v) &:= u_{n(v)}(v)\otimes e_{1 i_{n(v)}}
\end{align*} 
}with $I_v = \left\{i_2,\dots.i_{n(v)}\right\}$. Set $\widetilde{\mathcal{Z}}:=\widetilde{\mathcal{X}}\sqcup\left\{ y\otimes e_{11} : y \in \mathcal{Y}\right\}$ as a collection of elements of $\mathcal{G}\left(M^{n}\supseteq A^{n}\right)$ with $\widetilde{\mathcal{X}} := \left\{\tilde{u}_i(v) : i=1,\dots,n(v), v \in \mathcal{X}\right\}$, and 
\begin{equation*} 
P := 1\otimes e_{11} + \sum_{v\in\mathcal{X}}\sum_{i=1}^{n(v)-1} l\left(\tilde{u}_i(v)\right) = 1\otimes e_{11} + \sum_{v\in\mathcal{X}}\sum_{i=2}^{n(v)} r\left(\tilde{u}_i(v)\right).  
\end{equation*} 
By straightforward calculation we have  
\begin{equation}\label{eq2-Thm3-1} 
C_{\tau}\left(\mathcal{Z}\right)-1=C_{\tau^{n}}\big(\widetilde{\mathcal{Z}}\big)-\tau^{n}(P).  
\end{equation} 

(Step III: Reduction) Set $\widetilde{A} := A^{n}P$ and $\widetilde{M} := PM^{n}P$. Clearly, $\widetilde{M}$ is generated by the $g\otimes e_{11}$'s with $g \in \mathcal{G}$ and the $l\left(u_k(v)\right)\otimes e_{i_{k+1} 1}$'s with $k=1,\dots,n(v)-1$, $v \in \mathcal{V}$. Set $\tilde{\tau} := \tau^{n}\big|_{\widetilde{M}}$, and the $\tilde{\tau}$-conditional expectation $E_{\widetilde{A}}^{\widetilde{M}} : \widetilde{M} \rightarrow \widetilde{A}$ is given by the restriction of $E_A^M\otimes E_{\mathbf{C}^n}^{M_n(\mathbf{C})}$ to $\widetilde{M}$. Let $\widetilde{\mathcal{G}}$ be the smallest $E_{\widetilde{A}}^{\widetilde{M}}$-groupoid that contains $\big\{ g\otimes e_{11} : g \in \mathcal{G}\big\}\sqcup\big\{ l\left(u_k(v)\right)\otimes e_{i_{k+1} 1} : k=1,\dots,n(v)-1, v \in \mathcal{V}\big\}$. Also, for each $i=1,2,3$, let $\widetilde{\mathcal{G}}_i$ be the smallest $E_{\widetilde{A}}^{\widetilde{M}}$-groupoid that contains $\big\{ g\otimes e_{11} : g \in \mathcal{G}_i\big\}\sqcup\big\{ l\left(u_k(v)\right)\otimes e_{i_{k+1} 1} : k=1,\dots,n(v)-1, v \in \mathcal{V}\big\}$ and set $\widetilde{N}_i := \widetilde{\mathcal{G}}_i''$. Then, it is clear that   $\widetilde{N}_i = P\left(N_i\otimes M_n(\mathbf{C})\right)P = \left(N_i\otimes M_n(\mathbf{C})\right)\cap\widetilde{M}$, $i=1,2,3$. In particular, $\widetilde{A}$ is a Cartan subalgebra in $\widetilde{N}_3$, and thus $\widetilde{\mathcal{G}}_3 = \mathcal{G}\big(\widetilde{N}_3\supseteq\widetilde{A}\big)$. We have $\widetilde{\mathcal{G}} = \widetilde{\mathcal{G}}_1\vee\widetilde{\mathcal{G}}_2$, i.e., the smallest $E_{\widetilde{A}}^{\widetilde{M}}$-groupoid that contains both $\widetilde{\mathcal{G}}_1$ and $\widetilde{\mathcal{G}}_2$. Here, simple facts are in order.  
\begin{itemize}
\item[(a)] $N_1\otimes M_n(\mathbf{C})$ and $N_2\otimes M_n(\mathbf{C})$ are free with amalgamation over $N_3\otimes M_n(\mathbf{C})$ inside $M\otimes M_n(\mathbf{C})$ with subject to $E_{N_3}^M\otimes\mathrm{id}_{M_n(\mathbf{C})}$;  
\item[(b)] $v\otimes e_{11} = \left(a(v)^*\otimes e_{11}\right)\cdot \widetilde{u}_{n(v)}(v)\cdots\widetilde{u}_1(v)$ for every $v \in \mathcal{X}$; 
\item[(c)] $l\left(u_k(v)\right)\otimes e_{i_{k+1} 1} = \left(\widetilde{u}_k(v)\cdots \widetilde{u}_1(v)\right)\left(\left(u_k(v)\cdots u_1(v)\right)\otimes e_{11}\right)$ for every $k=1,\dots,n(v)-1$ and $v \in \mathcal{X}$. 
\end{itemize}
By (a), we have $\displaystyle{\widetilde{M} \cong \widetilde{N}_1\bigstar_{\widetilde{N}_3}\widetilde{N}_2}$ with respect to the restriction of $E_{N_3}^M\otimes\mathrm{id}_{M_n(\mathbf{C})}$ to $\widetilde{M}$ (giving the $\tilde{\tau}$-conditional expectation onto $\widetilde{N}_3$). By (b) and (c), it is plain to see that $\widetilde{\mathcal{Z}}$ is a graphing of $\widetilde{\mathcal{G}}$. Moreover, by its construction, it is decomposable, that is, $\widetilde{\mathcal{Z}} = \widetilde{\mathcal{Z}}_1\sqcup\widetilde{\mathcal{Z}}_2$ with collections $\widetilde{\mathcal{Z}}_1$, $\widetilde{\mathcal{Z}}_2$ of elements in $\widetilde{\mathcal{G}}_1$, $\widetilde{\mathcal{G}}_2$, respectively. Therefore, Lemma \ref{Lem3-11} shows that there is a treeing $\widetilde{\mathcal{Z}}' = \widetilde{\mathcal{Z}}'_1\sqcup\widetilde{\mathcal{Z}}'_2$ in $\widetilde{\mathcal{G}}_3 = \mathcal{G}\big(\widetilde{N}_3\supseteq\widetilde{A}\big)$ with the property that $\widetilde{\mathcal{Z}}_i\sqcup\widetilde{\mathcal{Z}}'_i$ is a graphing of $\widetilde{\mathcal{G}}_i$ for each $i=1,2$. Therefore, by Corollary \ref{Cor3-7} (b) (or Remark \ref{Rem3-8}), we have 
\begin{equation*}  
C_{\tilde{\tau}}\big(\widetilde{\mathcal{Z}}\big) 
= 
C_{\tilde{\tau}}\big(\widetilde{\mathcal{Z}}_1\sqcup\widetilde{\mathcal{Z}}'_1\big)+
C_{\tilde{\tau}}\big(\widetilde{\mathcal{Z}}_2\sqcup\widetilde{\mathcal{Z}}'_2\big)-
C_{\tilde{\tau}}\big(\widetilde{\mathcal{Z}}'\big) 
\geq 
C_{\tilde{\tau}}\big(\widetilde{\mathcal{G}}_1\big)+
C_{\tilde{\tau}}\big(\widetilde{\mathcal{G}}_2\big)-
C_{\tilde{\tau}}\big(\widetilde{\mathcal{G}}_3\big).  
\end{equation*} 
It is trivial that $c_{\widetilde{N}_i}(1\otimes e_{11}) = 1_{\widetilde{M}}$ for all $i=1,2,3$ with $1_{\widetilde{M}} = P$, and that $\left\{ g\otimes e_{11} : g \in \mathcal{G}_i\right\}=\left(1\otimes e_{11}\right)\widetilde{\mathcal{G}}_i\left(1\otimes e_{11}\right)$ for every $i=1,2,3$. Hence, \eqref{eq1-Thm3-1},\eqref{eq2-Thm3-1} and Propsotion \ref{Prop3-6} altogether imply that 
\allowdisplaybreaks{
\begin{align*} 
C_{\tau}(\mathcal{G}) + \varepsilon 
&\geq 
C_{\tau}(\mathcal{Z}) \\
&= 
C_{\tilde{\tau}}\big(\widetilde{\mathcal{Z}}\big) - \tilde{\tau}\big(1_{\widetilde{M}}\big)+1 \\
&\geq 
C_{\tilde{\tau}}\big(\widetilde{\mathcal{G}}_1\big)+
C_{\tilde{\tau}}\big(\widetilde{\mathcal{G}}_2\big)-
C_{\tilde{\tau}}\big(\widetilde{\mathcal{G}}_3\big)
- \tilde{\tau}\big(1_{\widetilde{M}}\big)+1 \\
&= 
C_{\tau}\big(\mathcal{G}_1\big)+
C_{\tau}\big(\mathcal{G}_2\big)-
C_{\tau}\big(\mathcal{G}_3\big). 
\end{align*}
}\end{proof}

\end{document}